\documentclass{article}
\usepackage{amssymb}
\usepackage{amsmath}
\usepackage{graphicx}   

\newtheorem{theorem}{Theorem}
\newtheorem{definition}{Definition}
\newtheorem{lemma}[theorem]{Lemma}
\newtheorem{rem}[theorem]{Remark}

\newcommand{\inte }{{\rm int}\,}

\newcommand{\mbeg}{\mathrm{b}}
\newcommand{\mend}{\mathrm{e}}

\newcommand{\Nbeg}{N^\mbeg}
\newcommand{\Nend}{N^\mend}

\newcommand{\Pibeg}{\Pi^\mbeg}
\newcommand{\Piend}{\Pi^\mend}
\newcommand{\rbeg}{r_\mbeg}
\newcommand{\rend}{r_\mend}
\newcommand{\barrbeg}{\bar{r}_\mbeg}
\newcommand{\barrend}{\bar{r}_\mend}

\newcommand{\xbeg}{x_\mbeg}
\newcommand{\xend}{x_\mend}
\newcommand{\ybeg}{y_\mbeg}
\newcommand{\yend}{y_\mend}
\newcommand{\tybeg}{\tilde{y}_u}

\newcommand{\stman}{x_s}
\newcommand{\unman}{y_u}

\newcommand{\pixend}{\pi_{\xend}}

\newcommand{\piyend}{\pi_{\yend}}

\newcommand{\dbeg}{d^\mbeg}
\newcommand{\dend}{d^\mend}

\newcommand{\Dbeg}{D^\mbeg}

\def\qed{{\hfill{\vrule height5pt width3pt depth0pt}\medskip}}

\title{\textbf{Connecting orbits for a singular nonautonomous real Ginzburg-Landau type equation}}
\author{{\Large Daniel Wilczak, Piotr Zgliczy\'nski\footnote{Research has been supported by Polish
National Science Centre grant 2011/03B/ST1/04780}}\\\\
\normalsize Faculty of Mathematics and Computer Science\\
\normalsize Jagiellonian University, \L ojasiewicza 6, 30-348 Krak\'ow\\
\normalsize \texttt{e-mail: \{Daniel.Wilczak, Piotr.Zgliczynski\}@ii.uj.edu.pl}
}
\begin{document}
\maketitle
\begin{abstract}
We propose a method for computation of stable and unstable sets associated to hyperbolic equilibria of nonautonomous ODEs and for computation of specific type of connecting orbits in nonautonomous singular ODEs. We apply the method to a certain a singular nonautonomous real Ginzburg-Landau type equation, which that arises from the problem of formation of spots in the Swift-Hohenberg equation.
\end{abstract}



\section{Introduction.}

The aim of this paper is to present a method for computation of connecting orbits between hyperbolic equilibria in nonautonomous ODEs. Although the algorithm we propose is general and may be applied to an ODE in any (finite) dimension, the main motivation for undertaking this research was the problem of finding transversal solutions for the following nonautonomous  real Ginzburg-Landau type equation
\begin{equation}
  A_{rr} + \frac{A_r}{r} - \frac{A}{4r^2}=A- A^3,  \label{eq:MS}
\end{equation}
where $A:[0,\infty) \to \mathbb{R}$ satisfies
\begin{equation}
A(0)=0, \quad \lim_{r \to \infty} A(r)=0. \label{eq:MS-asympt}
\end{equation}

This question  comes from McCalla and Sandstede work \cite{MS} on spots in the Swift-Hohenberg equation. The existence of solutions (\ref{eq:MS}) satisfying (\ref{eq:MS-asympt}) is listed there as an assumption in the construction of localized solutions called the B spots. The work of Scheel \cite{S} (references given there) indicates that bounded solutions of problem (\ref{eq:MS},\ref{eq:MS-asympt}) might be also interpreted as rotating wave solutions to reaction diffusion systems in the plane.

 It was pointed out in \cite{MS} that the  work by Scheel \cite{S} on (\ref{eq:MS},\ref{eq:MS-asympt}), where the existence of an infinite number of geometrically different solutions has been proved, contains a gap. Recently, a computer assisted proof of the existence of one such solution was given in \cite{BGW}.
The following theorem summarizes our main results about the system (\ref{eq:MS}).
\begin{theorem}\label{thm:six-orbits}
For $n=1,\ldots,6$ there exists a solution $A_n:(0,\infty)\to\mathbb R$ to (\ref{eq:MS},\ref{eq:MS-asympt}) such that:
\begin{itemize}
  \item $\lim_{r\to\infty} (A_n(r),A_n'(r)) = (0,0)$,
  \item $\lim_{r\to0^+} (A_n(r),A_n'(r)r) =  (0,0)$,
  \item The function $r\to A_n(r)$ has exactly $n$ local extremes in the domain $r\in(0,\infty)$.
\end{itemize}
The solution $A_n$ is transversal in the following sense.  The sets
\begin{eqnarray*}
W^u&=&\{ (r_0,A_0,A'_0)) \in (0,\infty) \times \mathbb{R} \times \mathbb{R} \ | \\
   & &\mbox{the backward solution of (\ref{eq:MS}) $A(r)$ with initial condition } \\
  & & \mbox{$A(r_0)=A_0$, $A'(r_0)=A'_0$ satisfies $\lim_{r \to 0} (A(r), A'(r)r)=0$}  \} \\
W^s&=&\{ (r_0,A_0,A'_0)) \in (0,\infty) \times \mathbb{R} \times \mathbb{R} \ | \\
  & &\mbox{the forward solution of (\ref{eq:MS}) $A(r)$ with initial condition }\\
 & & \mbox{$A(r_0)=A_0$, $A'(r_0)=A'_0$ satisfies $\lim_{r \to \infty} (A(r),A'(r))=0$}  \}
\end{eqnarray*}
are two dimensional immersed manifolds in the extended phase space $(r,A,A')$ and they intersect transversally along the solution curve defined by $A_n(r)$, i.e.
$(r,A_n(r),A_n'(r))$.
\end{theorem}

\begin{figure}
\centerline{\includegraphics[width=\textwidth]{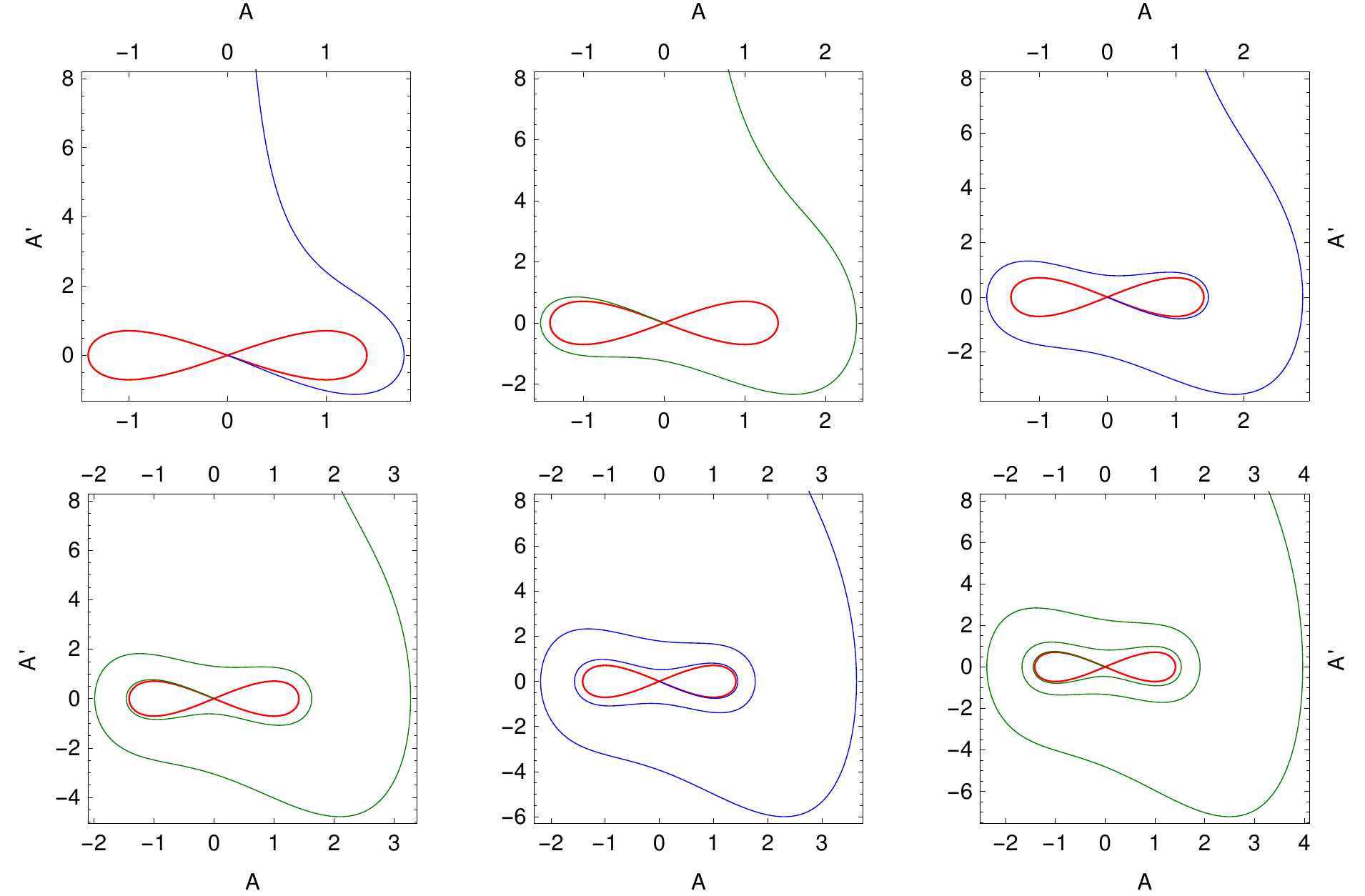}}
\centerline{\includegraphics[width=\textwidth]{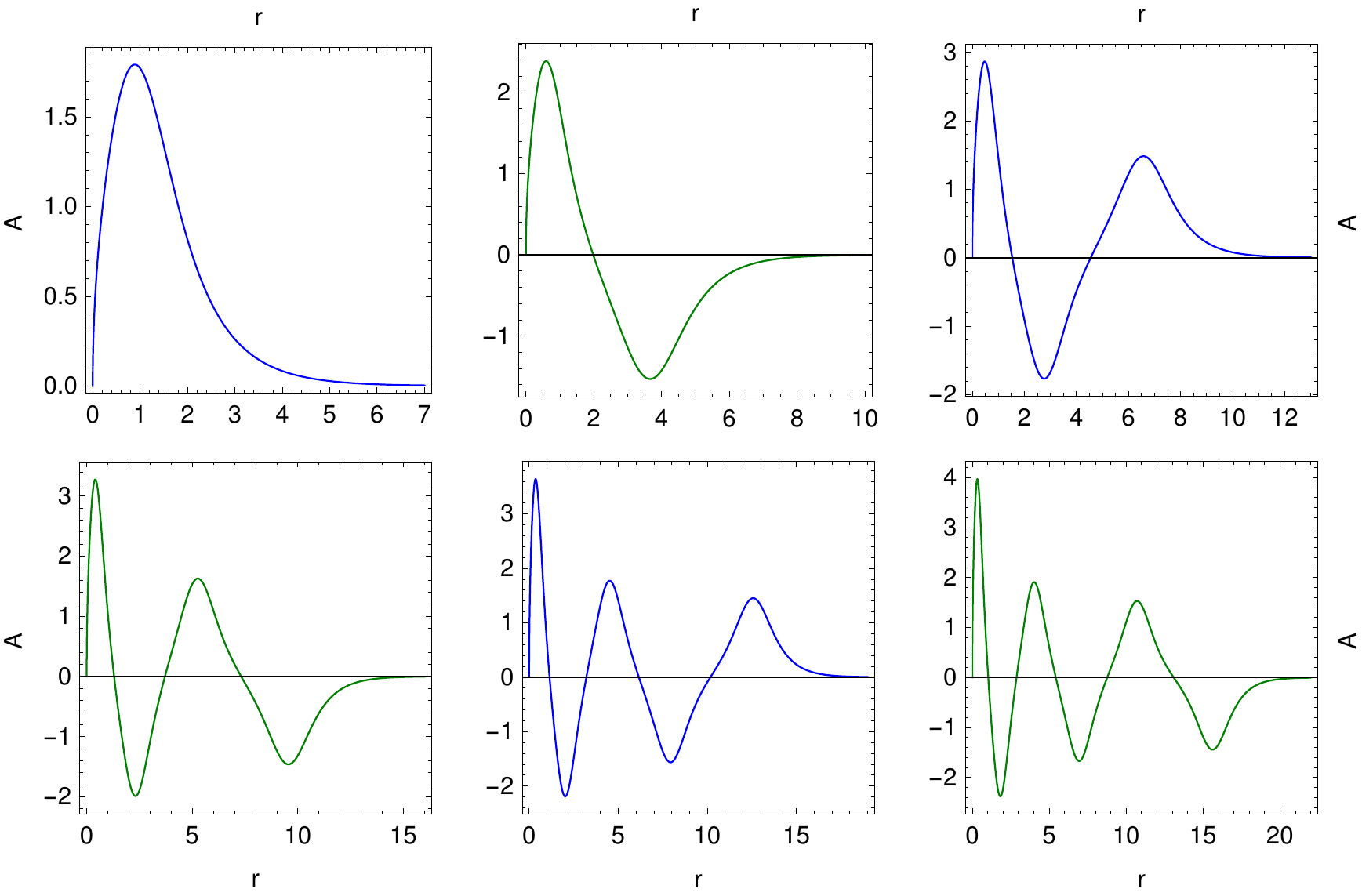}}
\caption{The six transverse connecting orbits resulting
from Theorem~\ref{thm:six-orbits}. In two upper rows in red we have plotted the separatrices (stable and unstable manifolds of $(0,0)$) of system $ A_{rr}=A- A^3$, which
is the formal limit of (\ref{eq:MS}) for $r \to \infty$.    \label{fig:six-orbits}}
\end{figure}

Plots of the six connecting orbits resulting from Theorem~\ref{thm:six-orbits} are shown in Fig.~\ref{fig:six-orbits}. The solutions differ in the number of local extremes or the number of sign changes. The orbit for $n=1$, with one maximum only, is the one already proved in \cite{BGW}. We present a proof of the existence of six such orbits, but  without  much difficulty it would be possible to validate the existence of several orbits with more  local extremes. The question of the existence of an infinite number of  geometrically different solutions for (\ref{eq:MS},\ref{eq:MS-asympt}) is still open.

It should be noted that Ryder \cite{R} proved the existence of  infinitely many geometrically different solutions to a similar problem
\begin{equation}
  u'' + \frac{2}{r}u' - u + u^3=0. \label{eq:ryder}
\end{equation}
plus boundary conditions analogous to (\ref{eq:MS-asympt}), which arises as the equation for the radial solutions in the Ginzburg-Landau equation in 3D. Apparently, Ryder did not proved the transversality.

The  Ryder's method applies to a wide class of equations, however due to the presence of  term $\frac{A}{4r^2}$, equation (\ref{eq:MS}) does not belong to this class.

Using our approach we can  obtain several first solutions also for (\ref{eq:ryder}) (apparently without much difficulty because it appears to be numerically easier than  (\ref{eq:MS})), but the number of solutions will be always finite in contrast with the Ryder's approach. Our shortcoming is due to the following:
the lack of mathematical understanding was replaced by the computer power. What is missing is the  understanding how the whole (un)stable set behaves away from the origin, and it is exactly the term $\frac{A}{4r^2}$, which causes the biggest problem.

Our approach is based on the shooting method (see for example \cite{HM}), which is a very well known technique in the dynamics of ODEs and boundary value problems for ODEs.

To realize the shooting method for problem (\ref{eq:MS},\ref{eq:MS-asympt}) we have to deal with the following issues
\begin{description}
\item[1.] the singularity for $r=0$,
\item[2.] the estimation of stable manifold of the point $(A,A')=(0,0)$ for non-autonomous equation (\ref{eq:MS}).
\end{description}
Issue 1 is dealt with via the change of the independent variable $\rho=\ln r$.
Then the equation (\ref{eq:MS}) becomes
\begin{equation}
  w" - w/4 = e^{2\rho}(w-w^3).  \label{eq:sing-r=0}
\end{equation}
Now we seek a solution such that $\lim_{\rho\to -\infty} w(\rho)=0$, so we need estimates for the unstable set of $(w=0,w'=0)$ for $\rho \to -\infty$ for our non-autonomous ODE (\ref{eq:sing-r=0}), which is basically issue 2 above.

Issue 2 in our work is split in two separate problems
\begin{itemize}
\item obtaining analytically bounds for the stable manifold in the neighborhood of the hyperbolic fixed point, which are valid up to some 'macroscopic' distance from it,
\item the globalization of the local (un)stable manifolds  using the rigorous numerical integrator from CAPD \cite{CAPD}.
\end{itemize}

Obtaining explicit bounds on the local (un)stable manifolds associated to hyperbolic equilibria in nonautonomous systems is an interesting issue itself and most likely this is the main mathematical contribution of this paper. For this end we adopt the method of cone conditions based on quadratic forms from \cite{ZCC,SZ} to the case of nonautonomous systems. The main novelty of this paper in this context is Lemma~\ref{lem:der-stable} which gives explicit bounds on the derivative of the parameterization of the stable manifold of a hyperbolic equilibrium with respect to the time variable. The method is general and may be applied to any smooth nonautonomous system in arbitrary finite dimension. We believe it should be an important ingredient when attacking the question of the existence of infinite number of geometrically different solutions to (\ref{eq:MS},\ref{eq:MS-asympt}) which is the next logical step for this project.

To demonstrate the effectiveness of our estimates of local (un)stable manifolds in Sections~\ref{sec:MSestm-stb-set} and ~\ref{sec:unstb-man-beg} we present easy calculations showing that using our approach we can estimate invariant manifolds for the system (\ref{eq:MS}) up to the distance of order $0.1$ from the origin.  Outside of this region the stable and unstable manifolds of (\ref{eq:MS}) and (\ref{eq:sing-r=0}) are continued to obtain a connecting orbit using a rigorous ODE integrator from the CAPD library \cite{CAPD}, but for this to succeed it is essential that our local bounds are valid up to a 'macroscopic' distance from the origin, otherwise the execution of the rigorous numerics will take very long time and the quality of bounds will be poor. Let us emphasize that to glue unstable and stable manifolds we need to integrate the system and associated variational equations over a finite time interval, only. Our local estimates on the invariant manifolds and bounds coming from the rigorous integration were sharp enough to fulfill the assumptions of the shooting method that guarantees the existence of six connecting orbits as stated in Theorem~\ref{thm:six-orbits}.

Let us comment about the computer assisted proof from \cite{BGW}. There the problem is formulated in the functional analytic language. Differential equation (\ref{eq:MS}) with the help of Green function is replaced by an integral equation on finite domain and conditions (\ref{eq:MS-asympt}) is encoded in the choice of the function space. Therefore obtaining a desired solution of (\ref{eq:MS},\ref{eq:MS-asympt}) is reduced to the fixed point problem for some nonlinear integral operator in some infinite dimensional function space. This should be contrasted with our approach, where we work in the phase space of our ODE.

The content of this paper can be described as follows. In Section~\ref{sec:Definitions-and-theorems} we recall the main geometric tools used in our approach to the (un)stable manifolds. In Section~\ref{sec:der-estm} we present theorems which provide us with bounds on the derivatives of (un)stable manifolds.
In Sections~\ref{sec:MSestm-stb-set}  and~\ref{sec:unstb-man-beg} we give explicit bounds for (un)stable manifolds and theirs derivatives for (\ref{eq:MS}) and (\ref{eq:sing-r=0}). In Section~\ref{sec:proof-for-MS} we give a computer assisted proof of Theorem~\ref{thm:six-orbits}.

\section{Geometric tools for invariant manifolds of fixed points}
\label{sec:Definitions-and-theorems}

To make this paper reasonably self-contained in this section we gather all necessary definitions and geometric theorems from \cite{ZGi,sym-per,ZCC,hom-tan} related to the rigorous investigation of invariant manifolds of fixed points, with modifications required for the setting of non-autonomous ODEs.

\subsection{Notation}
For a matrix $A \in \mathbb{R}^{n \times k}$ by $A^T$ we will denote its transpose. Let $Q$ be a square matrix. With some abuse of notation we will denote by the same letter a symmetric matrix $Q$ and associated quadratic form $Q(x)=x^T\cdot Q\cdot x$. To avoid ambiguity we will always use brackets $Q(\cdot)$ to indicate that $Q$ is treat as a quadratic form and $Qx$ or $Q\cdot (x_1-x_2)$ when we compute matrix-vector product. For a bilinear map $M:V \times W \to \mathbb{R}$ we will use the norm $\|W\|=\max_{\|v\|\leq 1, \|w\| \leq 1 } |M(v,w)|$.

For $u \in \mathbb{N}$, $p \in \mathbb{R}^u$ and $r>0$ by $B_u(p,r)$ we will denote an open ball of radius $r$ centered at  $p$. We will also often use $B_u=B_u(0,1)$. For a map $h:[0,1]\times X \to X$ and $t \in [0,1]$ we set $h_t(x)=h(t,x)$.

For a  function $f(t,x)$  we set $D_t f(t,x)=\frac{\partial f}{\partial t}(t,x)$ and $D_x(t,x)=\frac{\partial f}{\partial x}(t,x)$.

\subsection{Horizontal and vertical discs}

The main tools for proving the existence of local (un)stable manifolds in a given domain are isolating blocks and the cone conditions. The invariant manifolds are then constructed as graphs of smooth functions  in some isolating blocks, as  horizontal/vertical discs in h-sets. The size of an h-set gives an explicit range in which we can parameterize the local invariant manifold. The cone conditions are used to obtain both the existence and the local uniqueness of the invariant manifold as well as bounds on the derivatives of this parameterization.

\begin{definition} \cite[Def. 1]{ZCC}
\label{def:h-set}

An \emph{$h$-set} $N$ is a quadruple $\left(|N|,u(N),s(N),c_{N}\right)$ such that
\begin{itemize}
    \item $|N|$ is a compact subset of $\mathbb{R}^{n}$,
    \item $u(N),s(N)\in\{0,1,2,\ldots\}$ are such that $u(N)+s(N)=n$,
    \item $c_{N}:\mathbb{R}^{n}\rightarrow\mathbb{R}^{n}=\mathbb{R}^{u(N)}\times\mathbb{R}^{s(N)}$
    is a homeomorphism such that
    $c_{N}(|N|)=\overline{B}_{u(N)}\times\overline{B}_{s(N)}$.
\end{itemize}
We set

\[
\begin{array}{rcl}
    \mathrm{dim}(N)      & := & n,\\
    N_{c}       & := & \overline{B}_{u(N)}\times\overline{B}_{s(N)},\\
    N_{c}^{-}   & := & \partial B_{u(N)}\times\overline{B}_{s(N)},\\
    N_{c}^{+}   & := & \overline{B}_{u(N)}\times\partial B_{s(N)},\\
    N^{-}       & := & c_{N}^{-1}(N_{c}^{-}),\\
    N^{+}   & := & c_{N}^{-1}(N_{c}^{+}).
\end{array}
\]
\end{definition}

Hence an \emph{$h$-set} $N$ is a product of two closed balls in some coordinate system $c_N$.
 We call numbers $u(N)$ and $s(N)$ unstable and stable dimensions, respectively.
 The subscript $c$ refers to the new coordinates given by $c_{N}$. The set $|N|$ is called
 the support of an \emph{$h$-set}. We often drop bars in the symbol $|N|$ and use $N$
 to denote both the \emph{$h$-set} and its support.

Occasionally we will say that $N = \{(x_0,y_0)\} +
\overline{B}_u(0,r_1) \times \overline{B}_s(0,r_2) \subset
\mathbb{R}^u \times \mathbb{R}^s$ is an \emph{$h$-set}.  By this we will
understand a 'natural' \emph{$h$-set} structure on $N$ given by: $u(N)=u$,
$s(N)=s$, $c_{N}(x,y)= \left(\frac{x-x_0}{r_1},\frac{y-y_0}{r_2}
\right)$.  In the context of $\mathbb{R}^2$ and $u = 1$, $s = 1$ we
will sometimes write $N = z_0 + [-a,a] \times [-b,b]$.
This is compatible with the above
convention as $a$ defines the radius of the ball
$\overline{B}_u(0,a)=[-a,a]$ and $b$ of
$\overline{B}_s(0,b)=[-b,b]$.

\begin{definition} \cite[Def. 5]{ZCC}

Let $N$ be an \emph{$h$-set}. Let $b:\overline{B}_{u(N)}\rightarrow|N|$ be a continuous
mapping and let $b_{c}=c_{N}\circ b$. We say that $b$ is \emph{a horizontal disc} in $N$
if there exists a homotopy $h:[0,1]\times\overline{B}_{u(N)}\rightarrow N_{c}$, such that
\begin{eqnarray*}
     h_{0}      & = & b_{c}, \\
     h_{1}(x)   & = & (x,0),\quad\forall x\in\overline{u(N)}, \\
     h(t,x) \in N_{c}^{-}, && \forall t\in[0,1]\mathrm{\; and}\;\forall x\in\partial B_{u(N)}.
\end{eqnarray*}
\end{definition}

\begin{figure}
\centerline{\includegraphics[width=.4\textwidth]{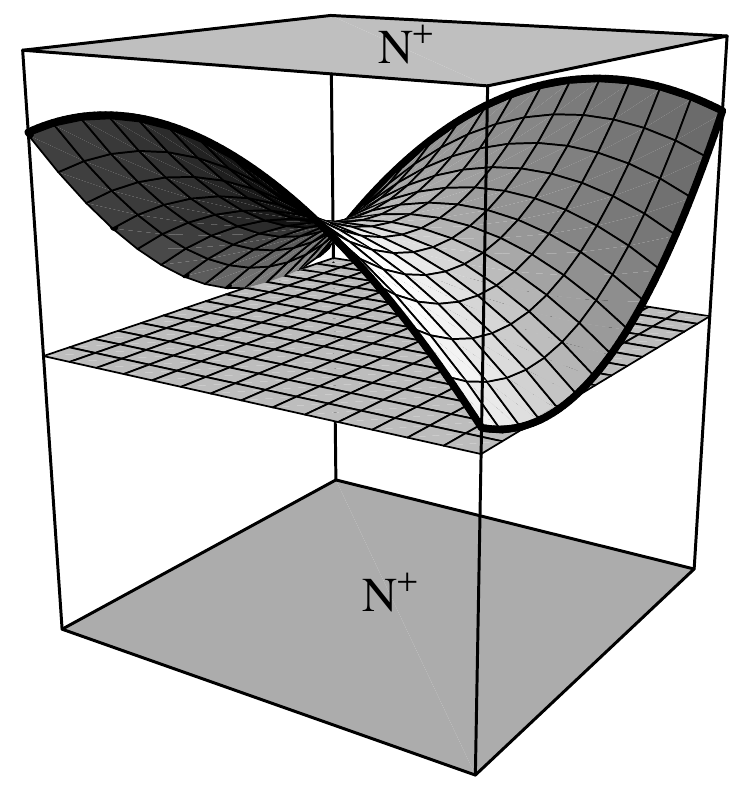}}
\caption{An h-set $N$ with $u=2$, $s=1$ and a horizontal disc in $N$.\label{fig:hdisc}}
\end{figure}

The geometry of this definition is shown in Fig.~\ref{fig:hdisc}. During the homotopy $h$ deforming a horizontal disc to the map $h(1,x)=(x,0)$
 the image of the boundary $h([0,1]\times B_{s(N)})$ must remain in $N_{c}^{-}$. In a similar way we can define vertical discs.

\begin{definition} \cite[Def. 6]{ZCC}
Let $N$ be an \emph{$h$-set}. Let $b:\overline{B}_{s(N)}\rightarrow|N|$ be a continuous
mapping and let $b_{c}=c_{N}\circ b$. We say that $b$ is \emph{a vertical disc}
in $N$ if there exists a homotopy $h:[0,1]\times\overline{B}_{s(N)}\rightarrow N_{c}$,
such that
\begin{eqnarray*}
    h_{0}       & = &   b_{c},\\
    h_{1}(x)    & = & (0,x),\quad\forall x\in\overline{B}_{s(N)},\\
    h(t,x)  \in N_{c}^{+}, && \forall t\in[0,1]\mathrm{\; and}\;\forall x\in\partial B_{s(N)}.
\end{eqnarray*}
\end{definition}

\begin{definition}\cite[Def. 7]{ZCC}
Let $N$ be an \emph{$h$-set} in $\mathbb{R}^{n}$and $b$ be a \emph{horizontal}
\emph{(vertical)} \emph{disc} in $N$. We will say that $x\in\mathbb{R}^{n}$ \emph{belongs}
to $b$, when $b(z)=x$ for some $z\in \mathrm{dom}(b)$.

By $|b|$ we will denote the image of $b$. Hence $z\in|b|\,$  iff
$z$ \emph{belongs to} $b$.
\end{definition}

\subsection{Cone conditions and the stable manifold theorem}

Below we present definitions and theorems that allow us to handle
and verify hyperbolic structures for nonautonomous ODEs using
$h$-sets and quadratic forms.

\begin{definition} \cite[Def. 8]{ZCC}
Let $N\subset\mathbb{R}^{n}$ be an \emph{$h$-set} and
$Q:\mathbb{R}^{n}\rightarrow\mathbb{R}$ be a quadratic form, such
that
\begin{equation*}
    Q(x,y) =
    \alpha(x)-\beta(y),\quad(x,y)\in\mathbb{R}^{u(N)}\times\mathbb{R}^{s(N)},
\end{equation*}
where $\alpha:\mathbb{R}^{u(N)}\rightarrow\mathbb{R}$ and
$\beta:\mathbb{R}^{s(N)}\rightarrow\mathbb{R}$ are positive
definite quadratic forms.

The pair $(N,Q)$ is called an \emph{$h$-set with cones}.
\end{definition}

We will often omit $Q$ in the symbol $(N,Q)$ and will say that $N$ is
an \emph{$h$-set with cones}.

\begin{definition} \cite[Def. 9]{ZCC}
\label{def:horizontal-disk-cone-condition} Let $(N,Q)$ be an
\emph{$h$-set with cones} and $b:\overline{B}_{u}\rightarrow|N|$
be a horizontal disc.

We will say that $b$ \emph{satisfies the cone condition (with
respect to $Q$)} iff for any $x_{1},x_{2}\in\overline{B}_{u}$,
$x_{1}\neq x_{2}$ the following inequality holds:
\begin{equation*}
    Q(b_{c}(x_{1})-b_{c}(x_{2})) > 0.
\end{equation*}
\end{definition}

\begin{figure}
\centerline{\includegraphics[width=.5\textwidth]{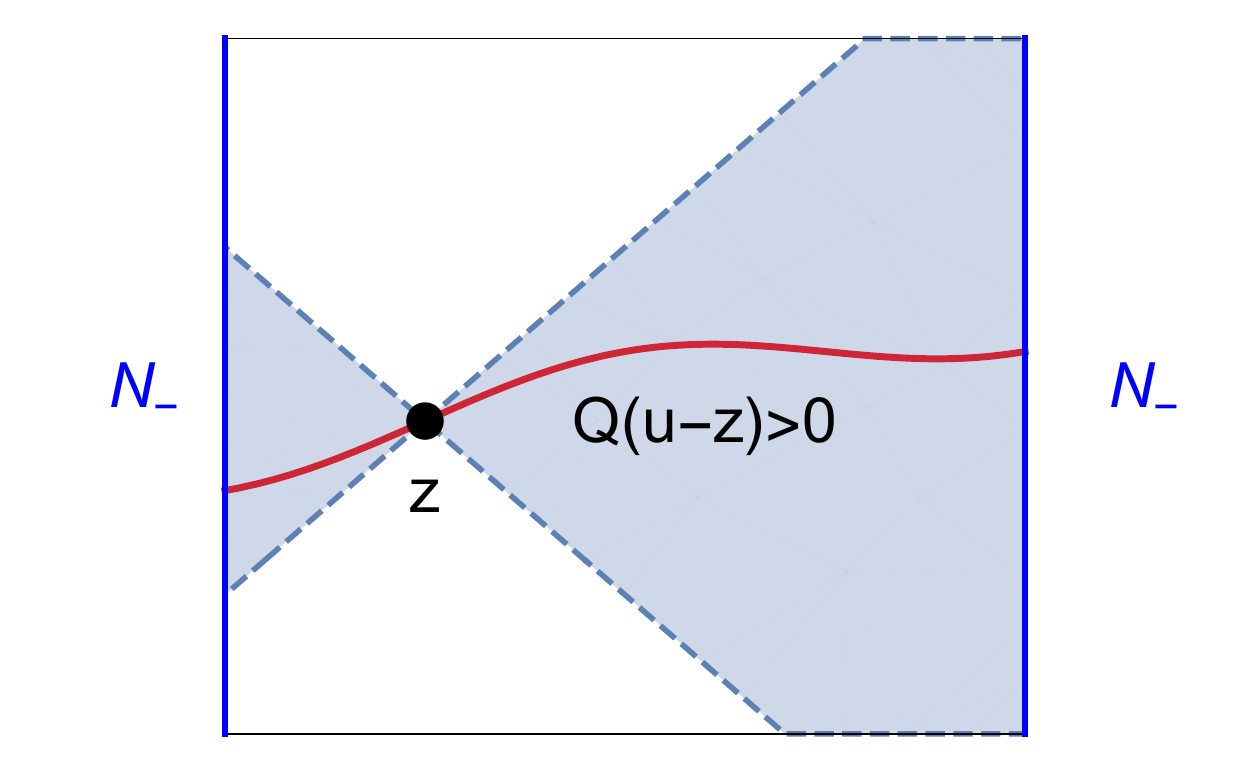}}
\caption{An h-set $N$ and a horizontal disc in $N$ satisfying the cone condition.\label{fig:hdiscWithCones}}
\end{figure}

The geometry of this definition is shown in Fig.~\ref{fig:hdiscWithCones}. The entire graph of  horizontal disc $b$ must lie in the positive cone attached to any point $z$ that belongs to the horizontal disc $b$ (except  point $z$).  It can be shown \cite{ZCC} that a horizontal disc satisfying the cone condition must be of the form $b_c(x,y)=(x,g(x))$ where $g$ is a Lipschitz function with an explicit bound on Lipschitz constant coming from the quadratic form. We will use this property to control derivative of the parameterization of the local invariant manifold which will be constructed as a horizontal (vertical) disc in some h-set satisfying the cone condition.

\begin{definition} \cite[Def. 10]{ZCC}
\label{def:vertical-disk-cone-condition}

Let $(N,Q)$ be an \emph{$h$-set with cones} and $b:\overline{B}_{s}\rightarrow|N|$ be
a vertical disc.

We will say that $b$ \emph{satisfies the cone condition (with
respect to $Q$)}  iff for any $y_{1},y_{2}\in\overline{B}_{s}$,
$y_{1}\neq y_{2}$ the following inequality holds:
\begin{equation*}
    Q(b_{c}(y_{1})-b_{c}(y_{2})) < 0.
\end{equation*}
\end{definition}

Let us consider an nonautonomous ODE
\begin{equation}
    \label{eq:non-auto}
    z' = f(t,z),
\end{equation}
where $z\in\mathbb{R}^{n}$, $f\in C^{1}(\mathbb{R} \times
\mathbb{R}^{n},\mathbb{R}^{n})$. Let us denote by $\varphi(t,t_0,p)$ the solution of
\eqref{eq:non-auto} with the initial condition $z(t_0)=p$.

Below we present the definition of the local stable and unstable set for a fixed point of an non-autonomous of ODE. These are just the
adaptation of the analogous notions for autonomous ODEs.

\begin{definition}
\label{def:hyperbolic-stable-unstable} Let
$z_{0}\in\mathbb{R}^{n}$ be such that $f(t,z_0)=0$ for $t \in
\mathbb{R}$. We will call such $z_0$ a fixed point for
(\ref{eq:non-auto}) (or $f$).

Let $Z\subset\mathbb{R}^{n}$, $z_{0}\in Z$ be a fixed point for
\eqref{eq:non-auto} and let $t_0 \in \mathbb{R}$. We define
\begin{eqnarray*}
    W_{t_0,Z}^s(z_{0}, \varphi)=W_{t_0,Z}^s(z_0,f)=
    \left\{ z:\forall_{t\ge 0}\varphi(t+t_0,t_0,z)\in Z,\ \lim_{t\rightarrow\infty}\varphi(t,t_0,z)=z_{0}\right\},\\
    W_{t_0,Z}^u(z_{0},\varphi)=W_{t_0,Z}^u(z_0,f)=
    \left\{ z:\forall_{t\leq 0}\varphi(t+t_0,t_0,z)\in Z,\ \lim_{t\rightarrow-\infty}\varphi(t,t_0,z)=z_{0}\right\}.
\end{eqnarray*}
When $\varphi$, $f$ or $z_0$ is known from the context, then we will
often drop it and write $W^{u,s}_{t_0,Z}(z_0)$ or $W^{u,s}_{t_0,Z}$.
\end{definition}

\begin{figure}
    \begin{center}
    \includegraphics[width=5.5cm]{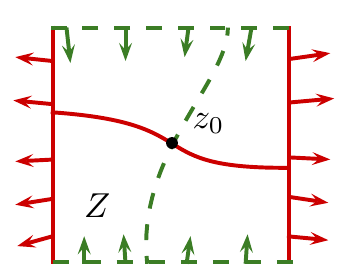}
    \end{center}
    \caption{\label{fig:stableunstable}An isolating block $Z$ for some planar ODE.
    The stable (vertical green dashed line) and unstable (horizontal red solid line)
    manifolds for $\varphi$ inside $Z$ are plotted, arrows indicate the vector field $f$.
   Dashed green and solid red border lines indicate the $\delta$-sections $\Sigma^+$ and $\Sigma^-$ respectively.  }
\end{figure}

Below we recall the notion of the isolating block from the Conley
index theory with some obvious modification for nonautonomous
ODEs.

\begin{definition}
\label{def:iso-block}  Let $\varphi$ be an autonomous  local flow
induced by $f:\mathbb{R}^n \to \mathbb{R}^n$. For $\delta>0$ the
set $\Sigma\subset\mathbb{R}^{n}$ is called a
$\delta\mathrm{-section}$ for the flow $\varphi$  iff
$\varphi\left((-\delta,\delta),\Sigma\right)$ is an open set and
the map $\sigma:\Sigma\times(-\delta,\delta)\rightarrow \varphi (
(-\delta,\delta),\Sigma)$ defined by $\sigma(x,t)=\varphi(t,x)$ is
a homeomorphism.

For a nonautonomous ODE $x'=f(t,x)$ inducing $\varphi$ we say that
$\Sigma$ as above is $\delta$-section for $\varphi$ for $t \in T
\subset \mathbb{R}$, iff for each $t_0 \in T$ the set $\Sigma$ is a
$\delta$ section for the vector field $x \rightarrow f(t_0,x)$ (i.e.
we freeze the time).

Let $B\subset\mathbb{R}^{n}$ be a compact set and $T \subset
\mathbb{R}$.  $B$ is called \emph{an isolating block} for
non-autonomous ODE for $t \in T$ iff $\partial B=B^- \cup B^+$,
where $B^-$ and $B^+$ are closed sets, and there exists $\delta>0$
and two $\delta\mathrm{-sections}$, $\Sigma^{+}$ and $\Sigma^{-}$
such that
\begin{eqnarray*}
    B^+ \subset \Sigma^+, \quad B^- \subset \Sigma^- , \\
   \forall x \in B^+ \quad  \forall t_0 \in T \quad \forall t \in (-\delta,0)  \qquad \varphi(t+t_0,t_0,x)
   \notin B, \\
    \forall x \in B^- \quad \forall t_0 \in T \quad  \forall t \in (0,\delta) \qquad \varphi(t,x)
   \notin B.
\end{eqnarray*}
$B^-$ and $B^+$ will  called an \emph{exit set} and an  \emph{entrance set}, respectively.
\end{definition}

In the present paper we will use \emph{$h$-sets} which are
isolating blocks for $t \geq t_0$ (or for $t \leq t_0$). Simply,
it means that $N^+$ and $N^-$ are sections of the vector field.

\begin{definition}
Let $N$ be an \emph{$h$-set} in $\mathbb{R}^n$. We say that $N$ is
an isolating block for nonautonomous ODE \eqref{eq:non-auto} for
$t \geq t_0$, iff $N^-$ and $N^+$ are $\delta$-sections for
$f(t,\cdot)$ for $t \geq t_0$ as in
Definition~\ref{def:iso-block}.
\end{definition}
Therefore, if  $N$, h-set, is an isolating block, then $N^-$ and $N^+$ are the exit and the entrance sets, respectively.

\begin{definition}
Let $N$ be an \emph{$h$-set}, such that $c_N$ is a diffeomorphism.
For a vector field $f$ on $|N|$ we define a vector field on $N_c$ by
\begin{equation}
  f_c(z)= Dc_N(c_N^{-1}(z))f(c_N^{-1}(z)).
\end{equation}
\end{definition}
Observe that $f_c$ is  the vector field $f$ expressed in
the new variables. If $f$ depends on time, then  $f_c$ is
time dependent too.

\begin{definition} \cite[Def. 13]{ZCC}

Let $U\subset\mathbb{R}^{n}$ be such that $U=\overline{U}$ and
$\inte{U} \ne \emptyset$. Let $g:U\rightarrow\mathbb{R}^{m}$ be a
$C^{1}$ function. We define the interval enclosure of $Dg(U)$ by
\begin{equation*}
    \left[Dg(U)\right]:=\left\{ A\in\mathbb{R}^{n\times m}:\forall_{i,j}A_{ij}\in\left[\inf_{x\in U}\frac{\partial g_{i}(x)}{\partial x_{j}},\sup_{x\in U}\frac{\partial g_{i}(x)}{\partial x_{j}}\right]\right\}.
\end{equation*}

We say that an interval matrix $[A]\subset \mathbb{R}^{n \times n}$ (i.e. a set of matrices)  is positive definite if for all symmetric matrices
$A\in[A]$  $A$ is positive definite.

\end{definition}

The lemma below serves  to motivate our definition of the cone condition for ODEs.
\begin{lemma}
\label{lem:cc-0DE}

Assume that $Q$ is a quadratic form on $\mathbb{R}^n$.

Consider equation (\ref{eq:non-auto}). Let $I \subset \mathbb{R}$ be an interval. Assume that $U \subset \mathbb{R}^n$ is convex and
the following condition is satisfied for $t \in I$
\begin{equation}
    \mathnormal{ matrix}\quad[D_zf(t,U)]^{T}Q+Q[D_z f(t,U)]\quad\mathnormal{is\ positive\
    definite.}  \label{eq:defodecc}
\end{equation}
Let $z_1(t)$ and $z_2(t)$ be two different solutions of (\ref{eq:non-auto}), such that for $t \in I$ $z_i(t) \in U$.

Then there exists $\eta >0$ such that
\begin{itemize}
\item[1.]
\begin{eqnarray}
  \frac{d}{dt} Q(z_1(t) - z_2(t)) &>& 0, \quad t \in I, \label{eq:ccOde-der-1} \\
  \frac{d}{dt} Q(z_1(t) - z_2(t)) &>& \pm \eta Q(z_1(t)-z_2(t)), \quad t \in I \label{eq:ccOde-der-Q}
\end{eqnarray}
\item[2.]
  if $t_0, t+t_0 \in I$, $t >0$ and $Q(z_1(t_0) - z_2(t_0)) >0$, then
    \begin{equation}
      Q(z_1(t_0+t) - z_2(t_0+t)) \geq e^{\eta t} Q(z_1(t_0)- z_2(t_0)), \label{eq:cc-exp-pos-cone}
    \end{equation}
\item[3.]
if $t_0, t+t_0 \in I$, $t <0$ and $Q(z_1(t_0) - z_2(t_0)) < 0$, then
    \begin{equation}
      Q(z_1(t_0+t) - z_2(t_0+t)) \leq e^{\eta |t|} Q(z_1(t_0)- z_2(t_0)). \label{eq:cc-exp-neg-cone}
    \end{equation}
\end{itemize}
\end{lemma}
\textbf{Proof:}
We have
\begin{eqnarray*}
 \frac{d}{dt}Q(z_1(t)-z_2(t))&=&
 (f(t,z_1(t)) - f(t,z_2(t)))^T\cdot Q\cdot (z_1(t)-z_2(t))  \\
  &+& (z_1(t)-z_2(t))^T\cdot  Q \cdot (f(t,z_1(t))-f(t,z_2(t))) \\
  &=&  (z_1(t)-z_2(t))^T \cdot C^T\cdot Q\cdot (z_1(t)-z_2(t)) \\
  &+& (z_1(t)-z_2(t))^T\cdot Q \cdot C(z_1(t)-z_2(t))\\
  &=&     (z_1(t)-z_2(t))^T\cdot (C^T Q + Q C)\cdot (z_1(t)-z_2(t)),
\end{eqnarray*}
where
\begin{eqnarray*}
C=C(t,z_1,z_2)=\int_0^1 D_z f(t,z_1 + s(z_2 -z_1))ds.
\end{eqnarray*}
Observe that $C \in [D_zf(t,U)]^{T}Q+Q[D_z f(t,U)]$,  so (\ref{eq:ccOde-der-1},\ref{eq:ccOde-der-Q}) follow easily.

From (\ref{eq:ccOde-der-Q}) and differential inequalities we obtain
\begin{eqnarray*}
  Q(z_1(t_0+t) - z_2(t_0+t)) &\geq& e^{\pm \eta t}  Q(z_1(t_0) - z_2(t_0)), \quad t>0, \ t_0,t_0+t \in I \\
  Q(z_1(t_0+t) - z_2(t_0+t)) &\leq& e^{\pm \eta t}  Q(z_1(t_0) - z_2(t_0)), \quad t<0, \ t_0,t_0+t \in I
\end{eqnarray*}
From the above conditions one can easily infer the remaining assertions.
\qed

Let us remark that from the above lemma it follows that  cone condition (\ref{eq:defodecc}) implies that cones $\mathcal{C}^+(z)=\{u \ | \ Q(z-u)>0\}$ are forward  invariant for (\ref{eq:non-auto}) and there is an exponential expansion in $\mathcal{C}^+(z)$. Analogously  cones $\mathcal{C}^-(z)=\{u \ | \ Q(z-u)>0\}$ are backward invariant and there is an exponential expansion in $\mathcal{C}^-(z)$ in backward time.

The following theorem can be easily obtained using the method as in the
proof of Theorem~26 in \cite{ZCC}.

\begin{theorem}
\label{thm:stb-man-thm}

Consider (\ref{eq:non-auto}).

Assume that $z_0 \in \mathbb{R}^n$ is such that  for any $t$ holds
$f(t,z_0)=0$.

Assume that $(N,Q)$ is an \emph{$h$-set with cones}, such that
$z_0 \in \inte N $ and which is an isolating block for
\eqref{eq:non-auto} for $t\geq t_0$, $c_N$ is a diffeomorphism and
that the following cone condition is satisfied for $t \geq t_0$:
\begin{equation}
    \mathnormal{ matrix}\quad[D_zf_c(t,N_c)]^{T}Q+Q[D_zf_c(t,N_c)]\quad\mathnormal{is\ positive\
    definite.}  \label{eq:odecc}
\end{equation}

Then for any $t_1\geq t_0$  the set $W_{t_1,N}^{s}(z_0,f)$ is a
vertical disc in $N$ satisfying the cone condition. Therefore there exists a function $x_s: \{t \geq t_0\} \times \overline{B}_s \to \overline{B}_u$, such that
\begin{equation}
  W_{t_1,N}^{s}(z_0,f)=\{c_N^{-1}(x_s(t_1,y),y)) \ | \ y \in \overline{B}_s\}.
\end{equation}
\end{theorem}
Analogous result is valid for the local unstable manifold.

Observe that from the standard results about exponential dichotomies (see for example \cite{C}) it follows that
function $x_s(t_1,y)$ is smooth with respect to $t_1$ and $y$ if $f$ and our coordinate change $c_N$
are smooth.

\section{The estimate for the derivative  with respect to time of
$W_{t,N}^{u,s}(z_0)$ for non-autonomous ODE}
\label{sec:der-estm}

In the sequel we  represent points in $\mathbb R^n$ by  $z=(x,y)$ with $x \in \mathbb{R}^u$ and $y \in
\mathbb{R}^s$,  $n=u+s$.

\subsection{The estimate for the time derivative  of $W_{t,N}^{u,s}(z_0)$}

The following lemma gives computable bounds on the first
derivatives of stable manifolds of a hyperbolic fixed point for a nonautonomous ODE with respect to the initial time. The proof is analogous to the proofs of the dependence on parameters  of the stable manifolds  of maps from  \cite[Thm. 21]{ZCC} with some refined estimates in \cite[Th. 4.1]{hom-tan} and for
parameterized ODEs \cite[Th. 5]{SZ}. The  theorem below and its proof is an adaptation of these results for nonautonomous ODEs.

\begin{lemma}
\label{lem:der-stable}

The same assumptions as in Theorem~\ref{thm:stb-man-thm}. $N$ has a natural h-set structure, $N=\bar{B}_u(p_x,r_x) \times \bar{B}_s(p_y,r_y)$ and that the quadratic form is given by (in natural coordinates) $Q(x,y)=a\|x\|^2 - b\|y\|^2$, where $a,b> 0$.

From the assumption on the cone condition for $Q$ there exists $E>0$ such that for $t \geq t_0$ and $z \in N$ there holds
\begin{equation}
  v^T \left(\left(\frac{\partial f}{\partial z}(t,z)\right)^T Q + Q\frac{\partial f}{\partial z}(t,z)\right) v \geq E \|v\|^2,
  \quad v \in \mathbb{R}^n. \label{eq:cc-non-auto}
\end{equation}

Put
\begin{equation}
  m=\sup_{t \geq t_0, z \in N}
  \left\|\left(\frac{\partial f}{\partial t}(t,z)\right)^T Q + Q \frac{\partial f}{\partial t}(t,z)\right\|
   \label{eq:norm-off-diag-terms}
\end{equation}
and
\begin{equation}
  \delta=\frac{a m^2}{E^2}. \label{eq:derWs-delta}
\end{equation}

Under above assumptions,
 if  $(x_i,y_i) \in W_{t_i,N}^s(z_0)$ and $t_i\geq t_0$ for $i=1,2$, then
\begin{equation}
  a\|x_1 - x_2\|^2 \leq \delta (t_1 - t_2)^2 + b\|y_1-y_2\|^2. \label{eq:derW-asser-negC}
\end{equation}

In particular, if $W_{t,N}^s(z_0)=\{ (\stman(t,y),y) \ | \ y \in
\overline{B}_s(p_y,r_y) \}$, then for $t_1,t_2 \geq t_0$ there holds
\begin{equation}
  \|\stman(t_1,y)-\stman(t_2,y)\| \leq \frac{m}{E} \, |t_1 - t_2|. \label{eq:derW-estm}
\end{equation}

\end{lemma}
\textbf{Proof:}
Consider the quadratic form on the extended phase-space $(t,x,y)$
\begin{equation*}
  \widetilde{Q}(t,z)=Q(z) - \delta t^2, \quad  \delta>0.
\end{equation*}
The value of $\delta$ will be fixed later during the proof (and it will turn out that (\ref{eq:derWs-delta}) is a good guess ).

Observe that, if for $i=1,2$ $t_i \geq t_0$, $\bar{z}_i \in W^s_{t_i,N}(f)$ and $z_i(t_i + t)$ for $t>0$ is a solution of (\ref{eq:non-auto}) such that $z_i(t_i)=\bar{z}_i$, then
\begin{equation}
  |\widetilde{Q}((t_1+t,z_1(t_1+t)) - (t_2+t,z_2(t_2+t))| \leq L, \quad \mbox{for $t \in \mathbb{R}_+$}. \label{eq:widet-Q-bounded}
\end{equation}
Indeed, since $z_i(t_i + t) \in N$ for $t >0$, so $|Q(z_1(t_1+t) - z_2(t_2+t))|$ is bounded, and $\delta ((t_1+t) - (t_2+t))^2=\delta (t_1-t_2)^2$ is constant.

We want find a positive $\delta$, such that the positive cone defined in terms of $\widetilde{Q}$ to be forward invariant for $t \geq t_0$ and we have an exponential growth
in this cone, i.e.
there exists $c >0$, such that
if $\bar{z}_1,\bar{z}_2 \in N$ and $t_1\geq t_0$, $t_2 \geq t_0$ are such that $\widetilde{Q}((t_1,\bar{z}_1) - (t_2,\bar{z}_2))>0$ and, $z_1(t+t_1)$ and $z_2(t+t_2)$ are solutions of (\ref{eq:non-auto})
 with $z_i(t_i)=\bar{z}_i$, $i=1,2$ , then for $t >0$  holds
 \begin{equation}
   \widetilde{Q}((t_1+t,z_1(t_1+t)) - (t_2+t,z_2(t_2+t)))>e^{c t} \widetilde{Q}((t_1,\bar{z}_1) - (t_2,\bar{z}_2)) \label{eq:der-estm-exp-growth}
 \end{equation}
 if  $z_1(t_1+s),z_2(t_2+s) \in N$ for $s \in [0,t]$.

This together with (\ref{eq:widet-Q-bounded}) implies that if $\bar{z}_i \in W^s_{t_i,N}(f)$ and $t_i \geq t_0$ for $i=1,2$, then
\begin{equation}
  \widetilde{Q}((t_1,\bar{z}_1) - (t_2,\bar{z}_2)) \leq 0.  \label{eq:wq-w-neg}
\end{equation}
From this condition and the computed expression for $\delta$ we will obtain our assertions.

So now we turn to finding $\delta$.

Let $z_1(t+t_1)$ and $z_2(t+t_2)$ be as in condition (\ref{eq:der-estm-exp-growth}). We set
\begin{multline*}
  L(t):=\widetilde{Q}((t+t_1,z_1(t+t_1)) - (t+t_2,z_1(t+t_2)))=\\
  Q(z_1(t+t_1)-z_2(t+t_2)) - \delta (t_1 - t_2)^2.
\end{multline*}
We have
\begin{multline*}
  L'(0)=
  \frac{d}{dt}  Q(z_1(t+t_1)- z_2(t+t_2))_{t=0}=\\
  (f(t_1,z_1(t_1))-f(t_2,z_2(t_2)))^T \cdot Q\cdot (z_1(t_1) - z_2(t_2)) + \\
  (z_1(t_1) - z_2(t_2))^T \cdot Q\cdot  (f(t_1,z_1(t_1))-f(t_2,z_2(t_2))).
\end{multline*}

Observe that
\begin{multline*}
  f(t_1,z_1(t_1))-f(t_2,z_2(t_2)) = \\
  \int_{0}^1 \frac{\partial f}{\partial t} ((1-s)(t_1,z_1(t_1)) + s(t_2,z_2(t_2)) )ds \cdot (t_1 - t_2) + \\
    \int_{0}^1 \frac{\partial f}{\partial z} ((1-s)(t_1,z_1(t_1)) + s(t_2,z_2(t_2)) )ds \cdot (z_1(t_1)-z_2(t_2)) = \\
    C_t (t_1-t_2) + C_z (z_1(t_1) - z_2(t_2)),
\end{multline*}
where
\begin{eqnarray*}
  C_t=C_t(t_1,z_1,t_2,z_2)= \int_{0}^1 \frac{\partial f}{\partial t} ((1-s)(t_1,z_1) + s(t_2,z_2)
  )ds,\\
  C_z=C_z(t_1,z_1,t_2,z_2)= \int_{0}^1 \frac{\partial f}{\partial z} ((1-s)(t_1,z_1) + s(t_2,z_2)
  )ds.
\end{eqnarray*}
From (\ref{eq:cc-non-auto}) we immediately obtain
\begin{equation*}
  v^T \left(C_z^T Q + QC_z\right) v \geq E \|v\|^2, \quad v \in
  \mathbb{R}^n.
\end{equation*}
From (\ref{eq:norm-off-diag-terms}) it follows that
\begin{equation*}
  m\geq \left\|C_t^T Q + Q C_t\right\|.
\end{equation*}
Continuing the previous derivation we obtain
\begin{multline*}
  L'(0) =
   (z_1-z_2)^T \left(C_z^T Q + Q C_z\right) (z_1 - z_2)  + \\
   (t_1-t_2)(C_t^T Q) (z_1-z_2) + (z_1-z_2)^T (QC_t) (t_1 - t_2).
\end{multline*}
Therefore we arrived at the following estimate
\begin{multline*}
  L'(0) = \frac{d}{dt} \widetilde{Q}((t+t_1,z_1(t+t_1)) - (t+t_2,z_2(t+t_2)))_{t=0} \geq\\ E \|z_1 - z_2\|^2 - m |t_1-t_2| \cdot \|z_1-z_2\|.
\end{multline*}
Now we make the use of the assumption $\widetilde{Q}((t_1,z_1)-(t_2,z_2)) \geq 0$. We have
\begin{eqnarray*}
  \sqrt{\frac{a}{\delta}} \,  \|z_1-z_2\| \geq  \sqrt{\frac{a}{\delta}} \, \|x_1-x_2\|  \geq |t_1 - t_2|.
\end{eqnarray*}
Therefore we obtain
\begin{eqnarray*}
   L'(0) \geq E \|z_1 - z_2\|^2 - m  \sqrt{\frac{a}{\delta}} \,  \|z_1-z_2\|^2 = \left( E-m  \sqrt{\frac{a}{\delta}} \right)  \|z_1-z_2\|^2.
\end{eqnarray*}
Hence the positive cone for $\widetilde{Q}$ will be forward invariant provided
\begin{equation*}
   E > m \sqrt{\frac{a}{\delta}},
\end{equation*}
which implies that
\begin{equation}
  \delta > \frac{a m^2}{E^2}. \label{eq:derW-delta}
\end{equation}
Let us fix any $\delta$ satisfying (\ref{eq:derW-delta}).

From the above reasoning we know that the positive cone for $\widetilde{Q}$ will be forward invariant and if $\widetilde{Q}((t_1,z_1)-(t_2,z_2)) \geq 0$, then
for any $t$, such that $z_i(t_i+t) \in N$ holds
\begin{eqnarray*}
  \frac{d}{dt}\widetilde{Q}((t+t_1,z_1(t+t_1))&-&(t+t_2,z_2(t+t_2))) \geq \left( E-m  \sqrt{\frac{a}{\delta}} \right)  \|z_1-z_2\|^2 \geq \\
    & &  c \widetilde Q((t_1+t,z_1(t_1+t))-(t_2+t,z_2(t_2+t))
\end{eqnarray*}
for some positive $c$. Therefore
\begin{equation*}
  \widetilde{Q}((t+t_1,z_1(t+t_1))-(t+t_2,z_2(t+t_2))) \geq \widetilde{Q}((t_1,z_1) - (t_2,z_2)) e^{ct}
\end{equation*}
for $t>0$ such that both $z_1(s+t_1)$ and $z_2(s+t_2)$ are in $N$ for $s \in [0,t]$.

Therefore we have established
(\ref{eq:wq-w-neg})  for any $\delta$ satisfying inequality (\ref{eq:derW-delta}), but passing to $\delta \to \frac{a m^2}{E^2}$ we
obtain  it also for $\delta = \frac{a m^2}{E^2}$.

If we set  $\bar{z}_i=(x_i,y_i)$, then from  (\ref{eq:wq-w-neg}) we have the following inequality 
\begin{equation*}
  a\|x_1 - x_2\|^2 \leq \delta (t_1 - t_2)^2 + b\|y_1-y_2\|^2.
\end{equation*}
This proves (\ref{eq:derW-asser-negC}). Inequality (\ref{eq:derW-estm}) is obtained, when we apply (\ref{eq:derW-asser-negC}) to $\bar{z}_i=(t_i,\stman(t_i,y),y)$.

\qed

In the application of Lemma~\ref{lem:der-stable} in this paper we will only use estimate (\ref{eq:derW-estm}). Moreover, since we are using the natural coordinates on $N$ cone condition (\ref{eq:odecc}) becomes
\begin{equation}
    \mathnormal{ matrix}\quad[D_z f(t,N)]^{T}Q+Q[D_z f(t,N)]\quad\mathnormal{is\ positive\
    definite.}  \label{eq:odecc-nat-coord}
\end{equation}
This simply means that we do not need to use the coordinate change $c_N$.

\subsection{Estimate for time derivative of $W_{t,N}^u(z_0)$}
By reversing the direction of time from Lemma~\ref{lem:der-stable} we obtain the following statement about the derivatives of the unstable manifold.
\begin{lemma}
\label{lem:der-unstable}
Assume that $f(t,z_0)=0$ for all $t\in\mathbb{R}$. Let $N$ has a natural h-set structure, $N=\bar{B}_u(p_x,r_x) \times \bar{B}_s(p_y,r_y)$,  $z_0 \in \inte N$ and that  quadratic form $Q$ on is given by (in natural coordinates) $Q(x,y)=a\|x\|^2 - b\|y\|^2$, where $a,b> 0$.

Assume that for $t\leq t_0$ the set $N$ is an isolating block for (\ref{eq:non-auto}) and $f$ satisfies  cone condition (\ref{eq:odecc}) on $N$ with
respect to quadratic form $Q$.

From the cone condition for $Q$ there exists $E>0$ such that for $t \leq t_0$ and $z \in N$
\begin{equation*}
  v^T \left(\left(\frac{\partial f}{\partial z}(t,z)\right)^T Q + Q\frac{\partial f}{\partial z}(t,z)\right) v \geq E \|v\|^2,
  \quad v \in \mathbb{R}^n.
\end{equation*}

Put
\begin{equation*}
  m=\sup_{t \leq t_0, z \in N}\left\|\left(\frac{\partial f}{\partial t}(t,z)\right)^T Q + Q \frac{\partial f}{\partial t}(t,z)\right\|. 
\end{equation*}

Under above assumptions,
 if $W_{t,N}^u(z_0)=\{ (x,\unman(t,x)) \ | \ x \in
\overline{B}_u(0,1) \}$, then for $t_1,t_2 \leq t_0$ holds
\begin{equation*}
  \|\unman(t_1,x)-\unman(t_2,x)\| \leq \frac{m}{E} \, |t_1 - t_2|.
\end{equation*}
\end{lemma}
\textbf{Proof:}
Transformation $t \to -t$ changes $Q \to -Q$ and $(x,y) \to (y,x)$.
The quantities in Lemma~\ref{lem:der-stable} are not affected by this change.
\qed

\section{Estimates for the stable manifold for $r \to \infty$ for equation (\ref{eq:MS})}
\label{sec:MSestm-stb-set}

In the limit $r \to \infty$ equation (\ref{eq:MS}) becomes
\begin{equation}
  u"=u-u^3.  \label{eq:end}
\end{equation}
We consider the first order system corresponding to (\ref{eq:end}).
\begin{equation}
  u'=v, \quad v'=u-u^3.  \label{eq:end-ord1}
\end{equation}
It is easy to see that (\ref{eq:end-ord1}) is Hamiltonian with $H(v,u)=\frac{v^2}{2} - \frac{u^2}{2} + \frac{u^4}{4}$.

It is clear that $(0,0)$ is a hyperbolic fixed for (\ref{eq:end-ord1}). Its stable and unstable manifolds
coincide to form an eight-shaped loop (called a separatrix) shown in two upper row drawings in Fig.~\ref{fig:six-orbits}.  Any solution of (\ref{eq:MS}) which converges to $0$ for $r \to \infty$ when plotted on $(A,A')$ plane  should approach this separatrix. In Fig.~\ref{fig:six-orbits} this is shown for six
orbits.

 We introduce coordinates, which diagonalize the linear part of (\ref{eq:end-ord1})
by
\begin{equation}
\label{eq:gc-n-out-inv}
\left\{\begin{aligned}
  \xend=u+v, \\
  \yend=u-v.
\end{aligned}
\right.
\end{equation}
The inverse is given by
\begin{equation*}
\left\{\begin{aligned}
  u=(\xend + \yend)/2, \\
  v=(\xend - \yend)/2.
\end{aligned}
\right.
\end{equation*}

In new coordinates system (\ref{eq:end-ord1}) becomes
\begin{equation}
\label{eq:end-diag}
\left\{\begin{aligned}
  \xend'=\xend + g(\xend,\yend), \\
  \yend'=-\yend - g(\xend,\yend)
\end{aligned}
\right.
\end{equation}
where
\begin{equation*}
  g(x,y)= - \frac{x - y}{2r} + \frac{x+y}{8r^2} - \frac{1}{8}(x + y)^3.
\end{equation*}


\subsection{Isolating block}
We look for the isolating block $\Nend=[-\dend_1,\dend_1]\times[-\dend_2,\dend_2]$ for
(\ref{eq:end-diag}).

\begin{lemma}\label{lem:isolation-end}
Put $d=\dend_1+\dend_2$ and assume that
\begin{equation}\label{eq:isolation-end}
  \min\{\dend_1,\dend_2\} > d\left(\frac{1}{2r}+\frac{1}{8r^2}+\frac{d^2}{8}\right).
\end{equation}
Then $\Nend$ is an isolating block for (\ref{eq:end-diag}).
\end{lemma}
\textbf{Proof:}

First we verify the exit condition. For $|\xend|=\dend_1$ we have
\begin{multline*}
  \xend' \xend = \xend^2 + \xend g(\xend,\yend)=\\
    \xend^2- \frac{1}{2} \frac{\xend-\yend}{r}\xend + \frac{(\xend+\yend)\xend}{8r^2} -
    \frac{(\xend+\yend)^3\xend}{8} \geq \\
    (\dend_1)^2- \frac{1}{2} \frac{\dend_1d}{r} - \frac{\dend_1 d}{8r^2} -
    \frac{d^3\dend_1}{8} = \\
    \dend_1\left(\dend_1- d\left(\frac{1}{2r} + \frac{1}{8r^2} +
    \frac{d^2}{8}\right)\right).
\end{multline*}
Therefore $\xend\xend'>0$ when (\ref{eq:isolation-end}) is satisfied.

In a similar way we verify the entry condition. For
$|\yend|=\dend_2$ we have
\begin{multline*}
  -(\yend' \yend) = \yend^2 + \yend g(\xend,\yend)\geq
    \dend_2\left(\dend_2- d\left(\frac{1}{2r} + \frac{1}{8r^2} +
    \frac{d^2}{8}\right)\right).
\end{multline*}
Therefore $-(\yend\yend)'>0$ when (\ref{eq:isolation-end}) is satisfied.
\qed

\subsection{Cone condition}
Our cones are defined in terms of the quadratic form
\begin{equation}
\label{eq:Qend}
  Q=\left[
      \begin{array}{cc}
        a & 0 \\
        0 & -1 \\
      \end{array}
    \right].
\end{equation}

Let
\begin{equation*}
  C(N,f)=[\overline{Df}^T(N)] Q + Q[\overline{Df}(N)].
\end{equation*}
where $[\overline{Df}(N)]$ is defined as follows: $M \in [\overline{Df}(N)]$ iff there exists a pair of points $z_1,z_2 \in N$ (can be equal) such that
$M=\int_0^1 Df(z_1+t(z_2-z_1))dt$.

\begin{lemma}\label{lem:conesOnConvexSets}
Assume that $N$ is a convex and compact set and $f:\mathbb{R}^n \to \mathbb{R}^n$  is $C^1$ smooth.
Assume that for every $z \in N$ there exists $\epsilon(z)\geq \epsilon_N >0$, such that for all $v \in \mathbb{R}^n$ there holds
\begin{equation*}
v^T(Df^T(z)Q+Q Df(z))v \geq \epsilon(z) \|v\|^2.
\end{equation*}

Then for every $M \in C(N,f)$ and $v \in \mathbb{R}^n$ there holds
\begin{equation*}
v^T(M^TQ+QM)v \geq \epsilon_N \|v\|^2.
\end{equation*}

\end{lemma}
\textbf{Proof:}
 For $z_1,z_2 \in N$ and $v\in \mathbb{R}^n$ we have
\begin{multline*}
  v^T \left(\left(\int_0^1  Df^T(z_1+t(z_2-z_1))tdt \right) Q + Q\left(\int_0^1  Df(z_1+t(z_2-z_1))dt \right) \right)v =\\
  \int_0^1 v^T ( Df^T(z_1+t(z_2-z_1))Q + Q  Df(z_1+t(z_2-z_1)))v dt \geq \epsilon_N \|v\|^2.
\end{multline*}
\qed

\begin{lemma}\label{lem:cone-end}
Let $\Nend=[-\dend_1,\dend_1]\times[-\dend_2,\dend_2]$ and $Q$ be given by
(\ref{eq:Qend}) with $a=1$. Assume that for some $E>0$ and $r_*>0$ there holds
\begin{equation*}
  1 > \frac{1}{2r_*} + \frac{3}{4}(\dend_1+\dend_2)^2  + E/2.
\end{equation*}
Then for all $v \in \mathbb{R}^n$, $M \in C(\Nend,f)$ and $r \geq r_*$
\begin{equation}\label{eq:cone-end}
v^T(MQ+QM)v \geq E \|v\|^2,
\end{equation}
i.e. the cone condition holds on $\Nend$ for $r \geq r_*$.
\end{lemma}
\textbf{Proof:}
According to Lemma~\ref{lem:conesOnConvexSets} it is sufficient to show that
the matrix $Df^T(\xend,\yend)Q+QDf(\xend,\yend)-E\cdot\mathrm{Id}$ is positive definite for each
$(\xend,\yend)\in\Nend$.
Put
\begin{eqnarray*}
  M_1=\left[
     \begin{array}{cc}
        1 & 0 \\
        0 & -1 \\
      \end{array}
  \right],
  \quad
 M_2= \left[
     \begin{array}{cc}
        1 & 1 \\
        -1 & -1 \\
      \end{array}
  \right],
  \quad
 M_3= \left[
     \begin{array}{cc}
        -1 & 1 \\
        1 & -1 \\
      \end{array}
  \right].
\end{eqnarray*}
We have
\begin{equation*}
  Df(\xend,\yend)=M_1 + \left(\frac{1}{8r^2} - \frac{3}{8} \left(\xend+\yend\right)^2\right)M_2 +
  \frac{1}{2r}M_3.
\end{equation*}

Easy calculations give
\begin{equation*}
  M_1^TQ=\left[
     \begin{array}{cc}
        a & 0 \\
        0 & 1 \\
      \end{array}
  \right],
  \quad
 M_2^TQ= \left[
     \begin{array}{cc}
        a & 1 \\
        a & 1 \\
      \end{array}
  \right],
  \quad
 M_3^TQ= \left[
     \begin{array}{cc}
        -a & -1 \\
        a & 1 \\
      \end{array}
  \right].
\end{equation*}
We obtain
\begin{multline*}
  Df^T(\xend,\yend) Q + Q Df(\xend,\yend) = \\ \left[
     \begin{array}{cc}
        2a & 0 \\
        0 & 2 \\
      \end{array}
  \right] + \left(\frac{1}{8r^2} - \frac{3}{8} \left(\xend+\yend\right)^2\right)
  \left[
     \begin{array}{cc}
        2a & 1+a \\
        a+1 & 2 \\
      \end{array}
  \right] + \\
  \frac{1}{2r}  \left[
     \begin{array}{cc}
        -2a & a-1 \\
        a-1 & 2 \\
      \end{array}
  \right].
\end{multline*}
Let us denote
\begin{equation}\label{eq:defSymbolW}
  w= \frac{1}{8r^2} - \frac{3}{8} (\xend+\yend)^2.
\end{equation}
Then
\begin{multline*}
   Df^T(\xend,\yend) Q + Q Df(\xend,\yend)-E\cdot\mathrm{Id}=\\\left[
   \begin{array}{cc}
     2a\left(1+w - \frac{1}{2r}\right)-E  & w(a+1) + \frac{(a-1)}{2r} \\
      w(a+1) + \frac{(a-1)}{2r} & 2\left(1+w+\frac{1}{2r} \right)-E
   \end{array}
   \right].
\end{multline*}
It is immediate that for the positive definiteness of this matrix we need
\begin{equation*}
  1+w - \frac{1}{2r}-E>0.
\end{equation*}

From now on we set $a=1$. The Ger\v{s}gorin Theorem gives us the following sufficient condition for the positive definiteness of the matrix $Df^T(\xend,\yend)Q+QDf(\xend,\yend)-E\cdot\mathrm{Id}$
\begin{eqnarray*}
  2\left( 1 + w - \frac{1}{2r}\right) -E > 2|w|,\\
   2\left( 1 + w + \frac{1}{2r}\right) - E > 2|w|.
\end{eqnarray*}
We end up with the following condition
\begin{equation*}
   1 + w - \frac{1}{2r} > |w|  + E/2.
\end{equation*}
Substituting back for $w$ defined by (\ref{eq:defSymbolW}) we obtain two inequalities that must be simultaneously satisfied
\begin{eqnarray*}
  1 + \frac{1}{8r^2} - \frac{3}{8}(\xend+\yend)^2 - \frac{1}{2r} &>&
  \frac{1}{8r^2} + \frac{3}{8}(\xend+\yend)^2 + E/2,\\
  1 - \frac{1}{2r} &>& \frac{3}{4}(\xend+\yend)^2  + E/2.
\end{eqnarray*}
Since for $(\xend,\yend) \in \Nend$ there holds $|\xend+\yend|\leq (\dend_1+\dend_2)$ we
obtain the following condition
\begin{equation*}
   1 > \frac{1}{2r} + \frac{3}{4}(\dend_1+\dend_2)^2  + E/2
\end{equation*}
and Lemma~\ref{lem:conesOnConvexSets} implies the assertion (\ref{eq:cone-end}).
\qed

\subsection{Estimate on  constant $m$}
\begin{lemma}\label{lem:m-estimate-end}
Let $\Nend= [-\dend_1,\dend_1]\times[-\dend_2,\dend_2]$ and put $d=\dend_1+\dend_2$. Then
for $z=(\xend,\yend)\in\Nend$ and $r \geq r_*$ there holds
  \begin{equation*}
    \left\|\left(\frac{\partial f}{\partial r}(r,z)\right)^T Q + Q \frac{\partial f}{\partial r}(r,z)\right\| \leq
     \frac{d}{2r_*^2}\left( 1 + \frac{3}{4 r_*} \right) \sqrt{a^2 + 1}.
  \end{equation*}
\end{lemma}
\textbf{Proof:}
Since we assumed that $Q$ is diagonal it is easy to see that
\begin{equation*}
   \left\|\left(\frac{\partial f}{\partial r}(r,z)\right)^T Q + Q \frac{\partial f}{\partial r}(r,z) \right\|= 2  \left\|Q \frac{\partial f}{\partial r}(r,z)\right\|.
\end{equation*}
We have
\begin{equation*}
  \frac{\partial f}{\partial r}(r,z)=
  \left( \frac{\xend-\yend}{2r^2} -
  \frac{3}{8r^3}(\xend+\yend)\right) \left[ \begin{array}{c} 1 \\ -1
    \end{array}
  \right],
\end{equation*}
hence
\begin{equation*}
  Q \frac{\partial f}{\partial r}(r,z) =
  \left( \frac{\xend-\yend}{2r^2} -
  \frac{3}{8r^3}(\xend+\yend)\right) \left[ \begin{array}{c} a \\ 1
    \end{array}
  \right].
\end{equation*}
Finally we obtain
\begin{eqnarray*}
   \left\|Q \frac{\partial f}{\partial r}(r,z)\right\|
   \leq \left( \frac{d}{2r_*^2} + \frac{3d}{8 r_*^3} \right) \sqrt{a^2 +
   1}=\frac{d}{2r_*^2}\left( 1 + \frac{3}{4 r_*} \right) \sqrt{a^2 + 1}.
\end{eqnarray*}

\qed

\subsection{Some numbers}
Consider  equation (\ref{eq:end-diag}) and quadratic form (\ref{eq:Qend}) with $a=1$ (this is the value used later in the actual proof). Lemma~\ref{lem:isolation-end} and Lemma~\ref{lem:cone-end} guarantee that the set $\Nend=[-\dend_1,\dend_1]\times [-\dend_2,\dend_2]$ is an isolating block for (\ref{eq:end-diag}) satisfying the cone condition for all $r\geq r_*$ provided the following inequalities hold true:
\begin{eqnarray*}
  \min\{\dend_1,\dend_2\} &>&  d\left(\frac{1}{2r_*}+\frac{1}{8r_*^2}+\frac{d^2}{8}\right),\\
 \frac{E}{2} + \frac{1}{2r_*} + \frac{3}{4}d^2 &<& 1,
\end{eqnarray*}
where $d=\dend_1+\dend_2$.

For example, if $r_* =2$ and $\dend_1=\dend_2=1/4$ then the first inequality holds true. From Lemma~\ref{lem:m-estimate-end} we get the bound $m\leq \frac{11}{128}\sqrt 2$.
Now the constant $E$ should satisfy $E<2-3/2d^2-1/r_*=9/8$. Take for instance $E=1$. This gives us a bound for the derivative
$|\stman'(r)|\leq \frac{m}{E}=\frac{11 \sqrt{2}}{128} \approx 0.121534$.

For $r$ large the second condition gives an upper bound for $d <
\frac{2}{3}\sqrt{3}$.

\section{Estimates for the unstable manifold for $\rho \to -\infty$ for equation (\ref{eq:MS})}
\label{sec:unstb-man-beg}
Let us rewrite here system (\ref{eq:sing-r=0})
\begin{equation}
  w" - w/4 = e^{2\rho}(w-w^3).  \label{eq:rho-neg}
\end{equation}
We want to study the unstable manifold of $(0,0)$ for $\rho < 0$.

Passing formally to the limit $\rho \to - \infty$ in
(\ref{eq:rho-neg}) we obtain the following linear system
\begin{equation}\label{eq:lin-rho-limit}
   w" - w/4 = 0.
\end{equation}

Let $p=w'$ and let us introduce new coordinates  $\xbeg,\ybeg$
\begin{equation}\label{eq:gc-n-in}
\left\{\begin{aligned}
  w=\xbeg+\ybeg, \\
  p=\frac{\xbeg-\ybeg}{2}.
\end{aligned}
\right.
\end{equation}
The inverse is given by
\begin{equation*}
\left\{\begin{aligned}
  \xbeg=\frac{w+2p}{2}, \\
  \ybeg=\frac{w-2p}{2}.
\end{aligned}
\right.
\end{equation*}
In the coordinates $(\xbeg,\ybeg)$ (\ref{eq:rho-neg}) becomes
\begin{equation}\label{eq:beg-ord1}
  \left\{\begin{aligned}
  \xbeg'=\frac{1}{2}\xbeg+e^{2\rho}g(\xbeg,\ybeg), \\
  \ybeg'=-\frac{1}{2}\ybeg - e^{2\rho}g(\xbeg,\ybeg),
\end{aligned}\right.
\end{equation}
where
\begin{equation*}
  g(\xbeg,\ybeg)=(\xbeg+\ybeg) - (\xbeg+\ybeg)^3.
\end{equation*}

\subsection{Isolating block}
\begin{lemma}
\label{lem:iso-blokc-beg}
Let $\Nbeg=[-d,d]^2$. If $d<1/2$ then $\Nbeg$ is an isolating block for (\ref{eq:beg-ord1}) for any $\rho$.
\end{lemma}
\textbf{Proof:} Observe that if $|\xbeg+\ybeg|\leq 1$ then $g(\xbeg,\ybeg)$ has
the same sign as $\xbeg+\ybeg$.  The sign of $\xbeg+\ybeg$ coincides with the sign
of the dominant coordinate. Therefore
\begin{eqnarray*}
\xbeg g(\xbeg,\ybeg)>0 & \text{for} & |\xbeg|=d, \\
\ybeg g(\xbeg,\ybeg)>0 & \text{for} & |\ybeg|=d.
\end{eqnarray*}

This implies that
\begin{eqnarray*}
  \xbeg'\xbeg = \frac{1}{2}\xbeg^2/2 +e^{2\rho}\xbeg g(\xbeg,\ybeg)\geq \frac{1}{2}\xbeg^2/2 & \text{if} & |\xbeg|=d,\\
  -\ybeg'\ybeg = \frac{1}{2}\ybeg^2/2 +e^{2\rho}\ybeg g(\xbeg,\ybeg)\geq \frac{1}{2}\ybeg^2/2 & \text{if} & |\ybeg|=d.
\end{eqnarray*}
\qed

\subsection{Cone condition}
In this section we will derive  inequalities that guarantee the cone condition on an isolating block for (\ref{eq:beg-ord1}).
Let us denote by $f$ the vector field (\ref{eq:beg-ord1}). We have
\begin{equation*}
 Df(\xbeg,\ybeg)=\left[\begin{array}{cc}
   \frac{1}{2} & 0 \\
   0 & -\frac{1}{2}
 \end{array}
 \right]
 + e^{2\rho}(1-3(\xbeg+\ybeg)^2) \left[\begin{array}{cc}
   1 & 1 \\
   -1 & -1
 \end{array}
 \right].
\end{equation*}
Put
\begin{equation}
\label{eq:Qbeg}
  Q=\left[
      \begin{array}{cc}
        a & 0 \\
        0 & -1 \\
      \end{array}
    \right].
\end{equation}

\begin{lemma}
Assume that the isolating block for (\ref{eq:beg-ord1}) has the form $\Nbeg=[-\dbeg_1,\dbeg_1]\times[-\dbeg_2,\dbeg_2]$.
\label{lem:cc-beg}
Put
\begin{equation*}
w=e^{2\rho}(1-3(\xbeg+\ybeg)^2).
\end{equation*}

If
\begin{eqnarray}
   1+2w &>& E > 0,       \label{eq:cc-beg-c1} \\
   a(1+2w)^2 - E (a+1)(1+2w) + E^2 &>& w^2(a+1)^2  \label{eq:cc-beg-c2}
\end{eqnarray}
for $(\xbeg,\ybeg)\in \Nbeg$ and $\rho\leq \rho_0$ then the cone
condition is satisfied on $\Nbeg$ for all $\rho\leq\rho_0$ for (\ref{eq:beg-ord1})
with the constant $E$.
\end{lemma}
\textbf{Proof:}
Denote by $f$ the vector field (\ref{eq:beg-ord1}). Direct computation gives
\begin{multline*}
  Df^T(\xbeg,\ybeg) Q + Q Df(\xbeg,\ybeg) = \\
  \left[
      \begin{array}{cc}
        a & 0 \\
        0 & 1 \\
      \end{array}
    \right] + e^{2\rho}(1-3(\xbeg+\ybeg)^2)
    \left[
      \begin{array}{cc}
        2a & a+1 \\
        a+1 & 2 \\
      \end{array}
    \right] = \\
 \left[
      \begin{array}{cc}
        a (1+2w) & w(a+1) \\
        w(a+1) & 1 + 2w \\
      \end{array}
    \right].
\end{multline*}
We will apply the standard criterion for positive definiteness of matrices to the matrix $M=( Df^T(\xbeg,\ybeg) Q + Q Df(\xbeg,\ybeg) - E \cdot \mathrm{Id} )$. We want $M_{22} > 0$ and $\det M >0$. The first condition gives (\ref{eq:cc-beg-c1}), the second reads
\begin{equation*}
  \left(a(1+2w)-E\right) \left(1+2w-E\right) - w^2 (a+1)^2 > 0
\end{equation*}
which is equivalent to (\ref{eq:cc-beg-c2}).
\qed

Using the Ger\v{s}gorin theorem \cite{G,V} we can derive conditions, which might more suitable for easy estimates to be done by hand. The positive definiteness of
$$
  \left[
      \begin{array}{cc}
        a (1+2w) & w(a+1) \\
        w(a+1) & 1 + 2w \\
      \end{array}
    \right]
$$
is implied by the following two inequalities
\begin{eqnarray}
  a(1+2w) - w(a+1) > E,  \label{eq:beg-a-w-gersh}\\
   1+2w - w(a+1) > E.\nonumber
\end{eqnarray}

This leads to the following lemma.
\begin{lemma}
\label{lem:cc-beg2}

Assume that $(\dbeg_1+\dbeg_2)^2 \leq \frac{1}{3}$,  $0<a \leq 1$ and the following condition is satisfied
\begin{eqnarray*}
 a(e^{2\rho_0}+1)-e^{2\rho_0} > E.
\end{eqnarray*}
Then the cone condition is satisfied on $\Nbeg$ for (\ref{eq:beg-ord1}) for $\rho \leq \rho_0$ with the constant $E$.
\end{lemma}
\textbf{Proof:}
Under the assumption on $\dbeg_1$ and $\dbeg_2$ we have
\begin{equation*}
w=e^{2\rho}(1-3(\xbeg+\ybeg)^2)\geq 0.
\end{equation*}
Since $0<a\leq 1$, the  positive definiteness of $Df(\xbeg,\ybeg)^TQ+QDf(\xbeg,\ybeg)$ is implied by only one condition
(\ref{eq:beg-a-w-gersh}) which reduces to
\begin{equation*}
  a(w+1) - w > E.
\end{equation*}
Observe that the range of $w$ for $(\xbeg,\ybeg) \in \Nbeg$ and $\rho \leq \rho_0$ is given by
$0< w\leq e^{2\rho_0}=w_0$. Since $a\leq1$, the function $w \mapsto a(1+w) - w$ is non-increasing, hence
it is enough to have this inequality for the largest possible value $w_0=e^{2\rho_0}$.

\qed

\subsection{Estimation of  constant $m$ when $\rho\to-\infty$.}
\begin{lemma}\label{lem:m-estimate-begin}
Denote by $f$ the vector field (\ref{eq:beg-ord1}). Assume that the isolating block for (\ref{eq:beg-ord1}) has the form $\Nbeg=[-\dbeg_1,\dbeg_1]\times[-\dbeg_2,\dbeg_2]$ with $d:=\max\{\dbeg_1,\dbeg_2\} \leq 1/2$. Assume the quadratic form is diagonal as in (\ref{eq:Qbeg}). Then for $z=(\xbeg,\ybeg) \in \Nbeg$ and $\rho \leq \rho_0$ there holds
\begin{equation*}
    \left\|\left(\frac{\partial f}{\partial \rho}(\rho,z)\right)^T Q + Q \frac{\partial f}{\partial \rho}(\rho,z)\right\| \leq 8d e^{2\rho_0} \sqrt{a^2+1}.
  \end{equation*}
\end{lemma}
\textbf{Proof:}
For $z=(\xbeg,\ybeg)\in\Nbeg$ we have
\begin{equation*}
  \frac{\partial f}{\partial \rho}(\rho,z)= 2 e^{2\rho} g(\xbeg,\ybeg) \left[ \begin{array}{c} 1 \\ -1
    \end{array}
  \right].
\end{equation*}
Since $Q$ is diagonal we get
\begin{equation*}
  Q \frac{\partial f}{\partial \rho}(\rho,z) = 2e^{2\rho}g(\xbeg,\ybeg)
  \left[ \begin{array}{c} a \\ 1  \end{array} \right].
\end{equation*}
Given that
\begin{equation*}
  |g(\xbeg,\ybeg) | = |\xbeg+\ybeg|(1-(\xbeg+\ybeg)^2) \leq |\xbeg+\ybeg| \leq 2d.
\end{equation*}
for $(\xbeg,\ybeg)\in\Nbeg$ we obtain
\begin{equation*}
  \left\|  Q \frac{\partial f}{\partial \rho}(\rho,z) \right\| \leq 4d e^{2\rho_0} \sqrt{a^2+1}.
\end{equation*}
\qed

\subsection{Some numbers}
We will show now that we can construct a quite large isolating block on which the cone condition is satisfied and still we can have good estimates on the derivative on parametrization of unstable set with respect to $\rho$.

From Lemma~\ref{lem:iso-blokc-beg} the set $\Nbeg=[-d,d]^2$ is an isolating block provided $d\leq 1/2$.

Lemma~\ref{lem:cc-beg2} guarantees the cone condition on $\Nbeg$ with  constant $E$ for $\rho\leq \rho_0$ provided the following two inequalities are satisfied
\begin{eqnarray*}
d &\leq& \frac{1}{2\sqrt3},\\
a(e^{2\rho_0}+1)-e^{2\rho_0} &>& E.
\end{eqnarray*}
In the case $a=1$ (this is the value used in our computer assisted proof) the last inequality reduces to $E<1$ and $\rho_0$ can be an arbitrary number. Finally, from Lemma~\ref{lem:m-estimate-begin} we obtain an estimate on the constant $m=8d e^{2\rho_0} \sqrt{2}$.

Explicit numbers satisfying all required inequalities:
\begin{itemize}
  \item Put $d =1/4$, $\rho_0=\ln 2$, $a=1$ and $E=1/2$. Then we have $m=8\sqrt{2}$, $|\unman'(\rho)| \leq m/E=16\sqrt2$.
  \item Put $d =1/4$, $\rho_0=\ln 2^{-4}$, $a=1$ and $E=1/\sqrt2$. Then we have $m=\sqrt{2}/128$, $|\unman'(\rho)| \leq m/E=1/64.$
\end{itemize}

\section{Proof of Theorem~\ref{thm:six-orbits}}
\label{sec:proof-for-MS}

\begin{figure}
\centerline{\includegraphics[width=.7\textwidth]{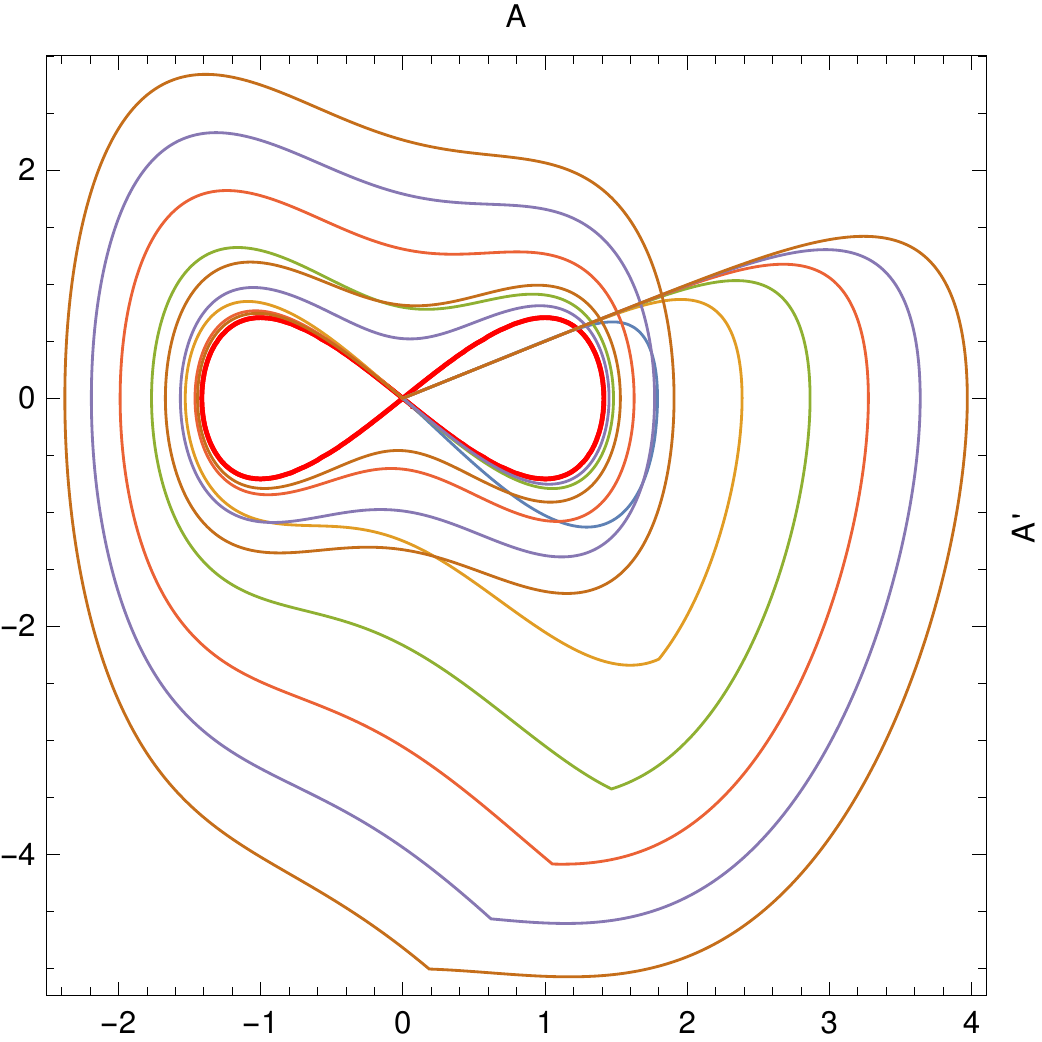}}
\caption{The six transverse connecting orbits resulting
from Theorem~\ref{thm:six-orbits}. The trajectories are plotted in two time
scales. For $r\geq1$ the points $(A_n(r),A_n'(r))$ are shown while for
$r\in(0,1)$ we plot $(A_n(\rho),A_n'(\rho))$ with $\rho=\ln r$. This figure
illustrates the convergence of orbits to $(0,0)$ when $\rho\to-\infty$. The figure-eight curve in the middle is the separatrix for (\ref{eq:end-ord1}).}\label{fig:six-orbits-scaled}
\end{figure}

In this section we give a computer assisted proof of the existence of six
geometrically different connecting orbits for equation (\ref{eq:MS}), i.e. solutions satisfying (\ref{eq:MS-asympt}).
Each orbit makes different number of (half) revolutions around the separatrix of equation (\ref{eq:end-ord1})
before approaching the equilibrium point $(A,A')=(0,0)$. Three of them converge
to $(0,0)$ from the right side ($A>0$) and three of them from the left side
($A<0$) - see Fig.~\ref{fig:six-orbits} and Fig.~\ref{fig:six-orbits-scaled}.

Before we give the proof of Theorem~\ref{thm:six-orbits} let us make several remarks about its content.
\begin{itemize}
\item $W^u$ can be expressed using other independent variable $\rho=\ln r$  and  $w(\rho)=A(e^\rho)$ (see Section~\ref{sec:unstb-man-beg}) as follows.
\begin{eqnarray*}
  W^u&=&\{ (\rho_0,w_0,w'_0)) \in \mathbb{R} \times \mathbb{R} \times \mathbb{R} \ | \\
   & &\mbox{such that the backward solution of (\ref{eq:rho-neg}) $w(\rho)$  with initial condition } \\
  & & \mbox{$w(\rho_0)=w_0$, $w'(\rho_0)=w'_0$ satisfies $\lim_{\rho \to -\infty} (w(\rho), w'(\rho))=0$}  \}.
\end{eqnarray*}
In fact in the proof we use the above description of $W^u$. The expression for $W^u$ given in the assertion of the theorem is obtained after we return
to the original variables $r$ and function $A(r)$.
\item The statement about the number of local extremes means that the trajectory $r\to(A_n(r),A_n'(r))$ intersects transversally the axis $A'=0$ exactly $n$ times.
\end{itemize}

In Fig.~\ref{fig:six-orbits-scaled} the six connecting orbits are shown in two
time scales: $r$ for $r\geq 1$ and $\rho=\ln r$ for $r\in(0,1)$. For $r=1$ they
coincide, i.e. $$(A(r=1),A'(r=1))=(A(\rho=0),A'(\rho=0))$$ therefore they appear
as continuous curves. This figure illustrates the convergence of orbits to
$(0,0)$ when $\rho\to-\infty$.

The proof of Theorem~\ref{thm:six-orbits} is based on the theorems introduced in the
previous sections. Using a non-rigorous simulation (the bisection) we found
good initial conditions for the connecting orbits. Then using validated numerics we verified the existence of isolating
blocks, cone conditions and conditions for the intersection of the  stable and unstable
manifolds. The details will be given in the next sections.

In Sections~\ref{sec:MSestm-stb-set} and~\ref{sec:unstb-man-beg}  we have found analytically the
isolating blocks on which the cone conditions are satisfied for the begin and end our connection. However, using
validated numerics we can find much tighter bounds for the
unstable/stable manifolds and similarly better bounds for the time-derivative of the parametrization of unstable/stable sets.

The input arguments to the algorithms that verify the existence of connecting
orbit are the following
\begin{itemize}
  \item $n$ - number of intersections of the trajectory with the line $A'=0$,
  \item $\hat r_n$ - an initial time,
  \item $\Delta r_n$ - time range around $\hat r_n$.
\end{itemize}
These parameters are listed in (\ref{eq:r0-data}) and found in nonrigorous
simulation (bisection).

\subsection{Isolating block and cone condition for $\rho\to-\infty$.}
\label{subsec:iso-bl-cc-beg}

By $(\xbeg,\ybeg)$ we denote coordinates that linearize the system at
$\rho\to-\infty$ as defined in (\ref{eq:gc-n-in}). We define an isolating block
\begin{equation}\label{eq:Nbeg_sizes}
\Nbeg = [-\dbeg_1,\dbeg_1]\times [-\dbeg_2,\dbeg_2] = [-0.125,0.125]\times[-2.8\cdot
10^{-8},2.8\cdot 10^{-8}].
\end{equation}
Let us define the following constants
\begin{equation}\label{eq:r0-data}
\begin{array}{|c|c|c|}
\hline n & \hat r_n & \Delta r_n\\ \hline
1 & 0.003288250 & 4\cdot 10^{-7}\\
2 & 0.001184020 & 6\cdot 10^{-8} \\
3 & 0.000650050 & 3\cdot 10^{-8}\\
4 & 0.000424204 & 2\cdot 10^{-8}\\
5 & 0.000304427 & 1\cdot 10^{-8}\\
6 & 0.000232050 & 8\cdot 10^{-9}\\ \hline
\end{array}
\end{equation}

\begin{lemma}\label{lem:Nbeg_estimates}
Put $\rho_*=\ln \left(\hat r_1+\Delta r_1\right)$, where
$\hat r_1$ and $\Delta r_1$ are defined in (\ref{eq:r0-data}).
\begin{itemize}
\item For $\rho\leq\rho_*$  set $\Nbeg$ is an isolating block for the
system (\ref{eq:beg-ord1}).
\item For $\rho\leq \rho_*$ the cone condition is satisfied on
$\Nbeg$ with quadratic form $Q(\xbeg,\ybeg)=\xbeg^2-\ybeg^2$ and thus for all
$\rho\leq\rho_*$   the set
\begin{equation*}
W^u_{\rho,\Nbeg}(0,0) = \{(\xbeg,\unman(\rho,\xbeg))\, |\,
\xbeg \in[-\dbeg_1,\dbeg_1]\}
\end{equation*}
 is a horizontal disc satisfying the cone condition.
\item The following estimate holds:
$$\left|\frac{\partial \unman}{\partial \rho}(\rho,\xbeg)\right|\leq
7.65\cdot 10^{-6}.$$
\end{itemize}
\end{lemma}

\textbf{Proof:}
Let $f(\rho,\xbeg,\ybeg)=(f_1(\rho,\xbeg,\ybeg),f_2(\rho,\xbeg,\ybeg))$ be the vector field as defined
in (\ref{eq:beg-ord1}) and put $f_\rho=f(\rho,\cdot,\cdot)$. Let $\dbeg_1,\dbeg_2$ be the sizes of $\Nbeg$ -- see
(\ref{eq:Nbeg_sizes}). Direct evaluation in interval
arithmetics gives the following inequalities
\begin{eqnarray*}
  f_1((-\infty,\rho_*]\times\{\dbeg_1\}\times[-\dbeg_2,\dbeg_2]) & \geq & 0.0625 > 0, \\
  f_1((-\infty,\rho_*]\times\{-\dbeg_1\}\times[-\dbeg_2,\dbeg_2])& \leq & -0.0625 < 0, \\
  f_2((-\infty,\rho_*]\times[-\dbeg_1,\dbeg_1]\times \{\dbeg_2\})& \leq &
  -4.8127929800192025\cdot 10^{-8} < 0, \\
  f_2((-\infty,\rho_*]\times[-\dbeg_1,\dbeg_1]\times \{-\dbeg_2\})& \geq &
  4.8127929800192025\cdot 10^{-8} > 0.
\end{eqnarray*}
This proves that $\Nbeg$ is an isolating block for all $\rho\leq \rho_*$.

Let us fix $\rho\in(-\infty,\rho_*]$ and recall that the quadratic form we are using is $Q(\xbeg,\ybeg)=\xbeg^2-\ybeg^2$.
Since $\dbeg_1+\dbeg_2\leq\frac{\sqrt{3}}{3}$ from Lemma~\ref{lem:cc-beg2}  we get that the cone condition is satisfied with any constant $E<1$.

Finally, from Lemma~\ref{lem:m-estimate-begin} and Lemma~\ref{lem:der-unstable} and passing to the limit with $E\to 1^-$
we get that
\begin{multline*}
\left|\frac{\partial \unman}{\partial \rho}(\rho,\xbeg)\right|\leq
	\frac{\left\|\frac{\partial f}{\partial \rho}((-\infty,\rho_*]\times \Nbeg)^TQ+Q\frac{\partial f}{\partial \rho}((-\infty,\rho_*]\times\Nbeg)\right\|}{E}\\
\leq \frac{4\sqrt2 e^{2\rho_*}|\dbeg_1+\dbeg_2|}{E}\quad \stackrel{E\to1^-}{\longrightarrow}\quad 4\sqrt2 e^{2\rho_*}|\dbeg_1+\dbeg_2|\\= 540761062255450812 \cdot10^{-23}\sqrt{2}< 7.65\cdot 10^{-6}.
\end{multline*}
\qed

\subsection{Isolating block and cone condition for $r\to\infty$.}

Let us consider an isolating block for $r\to\infty$. Put
\begin{equation}\label{eq:Nend_sizes}
\Nend = [-\dend_1,\dend_1]\times[-\dend_2,\dend_2] = [-0.0015,0.0015]\times[-0.01,0.01].
\end{equation}
\begin{lemma}\label{lem:Nend_estimates}
Let $r_*=6$.
\begin{itemize}
\item For $r\geq r_*$ the set $\Nend$ is an isolating block for the
system (\ref{eq:end-ord1}).
\item For $r\geq r_*$ the cone condition is satisfied on
$\Nend$ with quadratic form $Q(\xend,\yend)=\xend^2-\yend^2$ and thus for all
$r\geq r_*$ the set
\begin{equation*}
W^s_{r,\Nend}(0,0) = \{(\stman(r,\yend),\yend)\, |\,
\yend\in[-\dend_2,\dend_2]\}
\end{equation*}
 is a vertical disc satisfying the cone condition.
\item The following estimate holds:
$$\left|\frac{\partial \stman}{\partial r}(\rho,\yend)\right|\leq
0.000252.$$
\end{itemize}
\end{lemma}

\textbf{Proof:}
Let $f(r,\xend,\yend)=(f_1(r,\xend,\yend),f_2(r,\xend,\yend))$ be the vector field as
defined in (\ref{eq:end-ord1}) and let $\dend_1,\dend_2$ be the sizes of $\Nend$ --
see (\ref{eq:Nend_sizes}). Direct evaluation in interval
arithmetics gives the following inequalities
\begin{eqnarray*}
  f_1([r_*,\infty)\times\{\dend_1\}\times[-\dend_2,\dend_2]) & \geq &
  0.00051215277777777761 > 0,
  \\
  f_1([r_*,\infty)\times\{-\dend_1\}\times[-\dend_2,\dend_2])& \leq &
  -0.00051215277777777761 < 0, \\
  f_2([r_*,\infty)\times[-\dend_1,\dend_1]\times \{\dend_2\})& \leq &
  -0.010737706706597221 < 0, \\
  f_2([r_*,\infty)\times[-\dend_1,\dend_1]\times \{-\dend_2\})& \geq &
  0.010737706706597221 > 0.
\end{eqnarray*}
This proves that $\Nend$ is an isolating block for all $r\geq r_*$.

From the proof of Lemma~\ref{lem:cone-end} it follows that the cone condition
is satisfied on the set $\Nend$ for $r\geq r_*$
and with $a=1$, $E>0$ provided the following inequality holds true
$$
E<2-\frac{1}{r_*}-\frac{3}{2}(\dend_1+\dend_2)^2
$$
Substituting $r_*=6$ and $\dend_1,\dend_2$ form (\ref{eq:Nend_sizes}) we get
the upper bound
$$
E < E_0 = \frac{43995239}{24000000} \approx 1.8331349583333332.
$$

Finally, from Lemma~\ref{lem:m-estimate-end} and Lemma~\ref{lem:der-unstable}
we get that
\begin{multline*}
  \left|\frac{\partial \stman}{\partial r}(r,\yend)\right|\leq
  \frac{\left\|\frac{\partial f}{\partial r}([r_*,\infty)\times
  \Nend)^TQ+Q\frac{\partial f}{\partial r}([r_*,\infty)\times
  \Nend)\right\|}{E}\\ \leq
  \frac{\sqrt2(\dend_1+\dend_2)}{r_*^2}\left(1+\frac{3}{4r_*^2}\right)/E.
\end{multline*}
Passing to the limit with $E\to E_0^-$ we obtain
\begin{equation*}
  \left|\frac{\partial \stman}{\partial r}(r,\yend)\right|\leq
  \frac{\sqrt2(\dend_1+\dend_2)}{r_*^2}\left(1+\frac{3}{4r_*^2}\right)/E_0 = \frac{140875}{791914302}\sqrt2< 0.000252.
\end{equation*}
\qed

\subsection{Shooting between manifolds.}
Let us define three Poincar\'e sections. The first section is expressed in
coordinates $(r,\xbeg,\ybeg)$ -- see (\ref{eq:gc-n-in}) -- that linearize the system
at $\rho\to-\infty$. The section contains one of the exit edges of the set
$\Nbeg$, namely
\begin{equation*}
\Pibeg=\left\{(r,\xbeg,\ybeg) \,|\, \xbeg=\dbeg_1\right\}.
\end{equation*}
We will use $(r,\ybeg)$ coordinates to describe points on $\Pibeg$.

\begin{figure}
\centerline{
  \includegraphics[width=.7\textwidth]{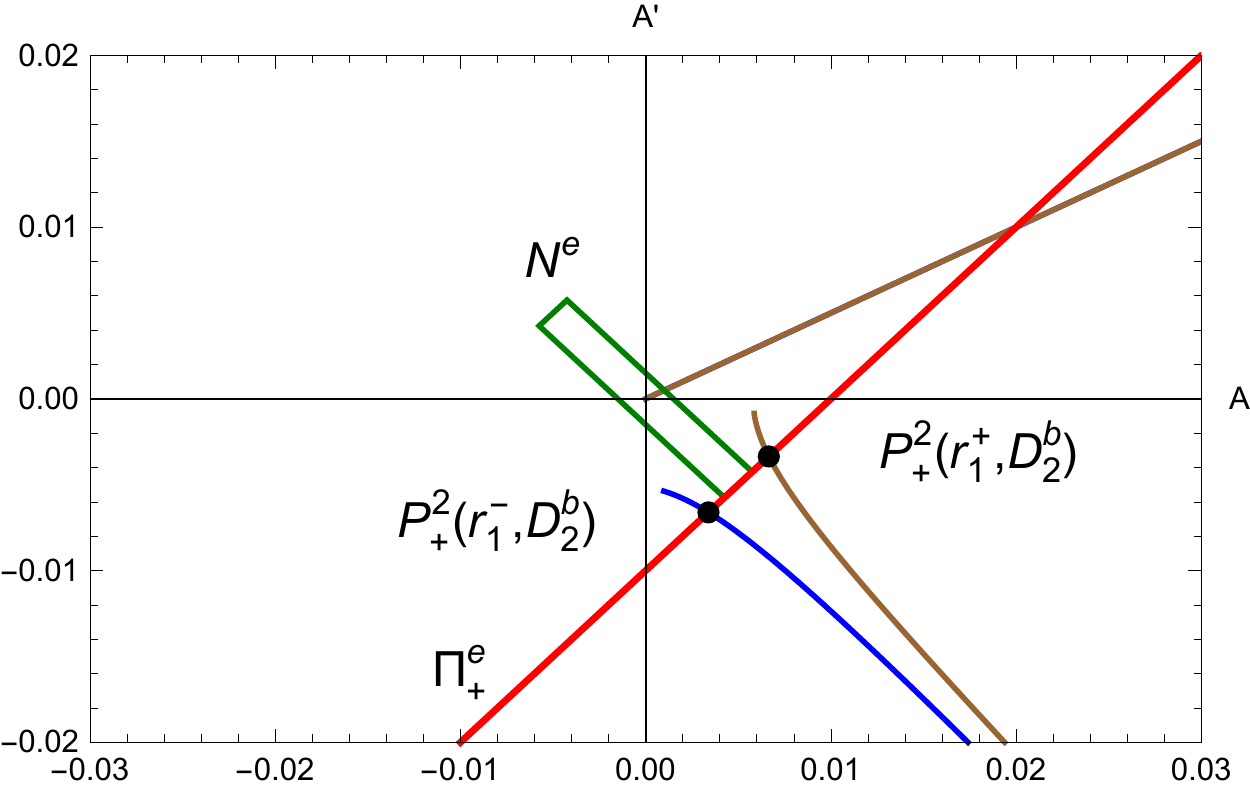}
}
\caption{The set $\Nend$, Poincar\'e section $\Piend_+$ and an indication that
$P^2_+(r_1^-,\Dbeg_2)$ and $P^2_+(r_1^+,\Dbeg_2)$ are mapped to opposite
sides of the set $\Nend$.\label{fig:shooting-first-orbit}}
\end{figure}

The two remaining sections are expressed in coordinates that define the set
$\Nend$ and each of them contains one of the entrance edges of this set
\begin{equation*}
\Piend_\pm=\left\{(r,\xend,\yend) \,|\, \yend=\pm \dend_2\right\}.
\end{equation*}
We will use $(r,\xend)$ coordinates to describe points on $\Piend_\pm$ when the sign
will be clear from the context. By $P_-:\Pibeg\to\Piend_-$ and
$P_+:\Pibeg\to\Piend_+$ we denote two Poincar\'e maps.
By $P^n_{\pm}(r,\ybeg)$ we denote the $n$-th itersection of the trajectory of (\ref{eq:MS}) starting at $(r,\dbeg_1,\ybeg)$ with the section $\Piend_\pm$. In what follows we will always
use odd number of intersections for the mapping $P^n_-$ and even for the mapping
$P^n_+$. Therefore we will skip the subscript $\pm$ to simplify the notation and
write just $P^n$.
\begin{lemma}\label{lem:shooting}
Let $\hat r_n$, $\Delta r_n$, $n=1,\ldots,6$ be as defined in
(\ref{eq:r0-data}) and put
\begin{eqnarray*}
r^-_n&=&\hat r_n -\Delta r_n\\
r^+_n&=&\hat r_n +\Delta r_n\\
\Dbeg_2 &=& [-\dbeg_2,\dbeg_2].
\end{eqnarray*}
Then for $n=1,\ldots,6$ the mapping
$P^{n+1}$ is well defined and smooth on the set
$$[r^-_n,r^+_n]\times [-\dbeg_2,\dbeg_2].$$
Moreover
\begin{eqnarray*}
  P^{n+1}(\{r^\pm_n\}\times \Dbeg_2) &\subset& \left\{(r,\xend)\,:\, |\xend|
  > \dend_1, r\geq r_*=6\right\},\\
\pi_x P^{n+1}(r^-_n,\ybeg)\cdot \pi_x P^{n+1}(r^+_n,\ybeg)&<&0\quad\text{for }
\ybeg\in \Dbeg_2
\end{eqnarray*}
and each trajectory of a point $(r_\mbeg,\ybeg)\in[r^-_n,r^+_n]\times
[-\dbeg_2,\dbeg_2]$ intersects the axis $A'=0$ exactly $n$ times before reaching the
section $\Piend_{(-1)^n}$.
\end{lemma}
\textbf{Proof:}
The proof is computer assisted and  uses an rigorous ODE solver from the CAPD
library \cite{CAPD}. In (\ref{eq:return-time}) we give  rigorous bounds for the
return time, i.e. for all $(r_\mbeg,\ybeg)\in [r^-_n,r^+_n]\times [-\dbeg_2,\dbeg_2]$ the image
$(r_\mend,\xend):=P^{(n+1)}(r_\mbeg,\ybeg)$ exists and $r_\mend$ belongs to the corresponding
interval listed in (\ref{eq:return-time}). This shows that such an intersection occurs always for $r_\mend\geq
r_*=6$.
\begin{equation}\label{eq:return-time}
\begin{array}{|c|c|}
\hline n & \text{bound for the return time }\rend \\ \hline
1 & [6.5694270711914049, 6.8663028711914071]\\
2 & [9.547364685097655, 9.8754898050976578] \\
3 & [12.63188037430908, 12.975630434309084]\\
4 & [15.467339314615232, 15.811089354615238]\\
5 & [18.659449295387446, 19.009330830371958]\\
6 & [21.645791630860828, 21.989541646860836]\\ \hline
\end{array}
\end{equation}
In (\ref{eq:cover}) we list estimates which establish the remaining
inequalities. We see that the $x$-coordinate of
$P^{n+1}(\{r_n^\pm\}\times \Dbeg_2)$ has opposite sign for
$r_n^-$ and $r_n^+$ and in each case its absolute value is bigger than
$\dend_1=0.0015$.
\begin{equation}\label{eq:cover}
\begin{array}{|c|c|c|}
\hline n &  \pi_xP^{n+1}(\{r_n^-\}\times \Dbeg_2)
& \pi_xP^{n+1}(\{r_n^+\}\times \Dbeg_2)
 \\
\hline
1  &  -0.0032[19,26] &  0.0032[89,97]\\
2  & 0.00189[2,8]  & -0.0018[86,93]\\
3  & -0.0019[87,92]  & 0.0019[37,43]\\
4  & 0.00209[5,9]  & -0.00221[0,5]\\
5  & -0.00155[5,9]  & 0.0015[78,82]\\
6  & 0.00173[1,4]  & -0.0016[57,61]\\ \hline
\end{array}
\end{equation}
In Fig.~\ref{fig:shooting-first-orbit} we show the isolating block $\Nend$,
Poincar\'e section $\Piend_+$ and by black dots we marked the image
$P^2_+(\{r_1^\pm\}\times \Dbeg_2)$.
\qed

\textbf{Proof of Theorem~\ref{thm:six-orbits}:}

Recall that we used $(\rbeg,\ybeg)$ coordinates on section $\Pibeg$ and $(\rend,\xend)$
coordinates on $\Piend_+$. From Lemmas~\ref{lem:Nbeg_estimates} and ~\ref{lem:Nend_estimates} we know that
\begin{eqnarray*}
W^u_{\rho,\Nbeg}(0,0) = \{(\xbeg,\unman(\rho,\xbeg))\, |\, \xbeg \in[-\dbeg_1,\dbeg_1]\}, \\
W^s_{r,\Nend}(0,0) = \{(\stman(r,\yend),\yend)\, |\,\yend\in[-\dend_2,\dend_2]\},
\end{eqnarray*}
where $y_u(\cdot,\cdot)$, $x_s(\cdot,\cdot)$  are  smooth functions. Put $$\tybeg(\rbeg,\xbeg):=\unman(\ln \rbeg,\xbeg)$$ and
consider the mapping
\begin{equation*}
  F_n(\rbeg)=\pixend P^{n+1}\left(\rbeg,\tybeg(\rbeg,\dbeg_1)\right)
              -x_s\left(\pi_{r} P^{n+1}(\rbeg,\tybeg(\rbeg,\dbeg_1)),\dend_2\right).
\end{equation*}

Zeros of this function correspond to connecting orbits we
are searching for. We will show that for $n=1,\ldots,6$ this function has a
unique  zero in $[r_n^-,r_n^+]$.

From Lemma~\ref{lem:shooting} the function $F_n$ is continuous and $F_n(r_n^-)F_n(r_n^+)<0$, because $|x_s\left(\pi_{r} P^{n+1}(\rbeg,\tybeg(\rbeg,\dbeg_1)),\dend_2\right)| \leq d_1^e$.
This proves that $F_n$ has zero in $[r_n^-,r_n^+]$. For uniqueness it is enough to show that $F_n'(r)\neq 0$ for
$r\in[r_n^-,r_n^+]$.

The proof is computer assisted and we will give details for $n=1$, only. For $n=2,\ldots,6$ we will
give the necessary estimates. Using the CAPD library \cite{CAPD} we computed a rigorous bound for the
derivative of the return map $P_+^2$ on the set $W=[r_1^-,r_1^+]\times[-\dbeg_2,\dbeg_2]$ and we got
\begin{eqnarray*}
DP^2_+(W) &=&
\begin{bmatrix}
\frac{\partial \pi_{r} P^2_+}{\partial \rbeg}(W) & \frac{\partial \pi_{r} P^2_+}{\partial \ybeg}(W)\\
\frac{\partial \pixend P^2_+}{\partial \rbeg}(W) & \frac{\partial \pixend P^2_+}{\partial \ybeg}(W)
\end{bmatrix}\\
&\subset&
\begin{bmatrix}
[-292366, 565586] & [-58.3544, 88.8462]\\
[1538.59, 16826.7] &[-0.253485, 2.36953]
\end{bmatrix}.
\end{eqnarray*}

From Lemma~\ref{lem:Nbeg_estimates} and Lemma~\ref{lem:Nend_estimates} we have
bounds for the partial derivatives of the parametrization of invariant sets
with respect to time variable, namely
\begin{eqnarray*}
\left|\frac{\partial \unman}{\partial \rho}(\rho,\dbeg_1)\right|&\leq&
7.65\cdot 10^{-6},\\
\left|\frac{\partial \stman}{\partial r_\mend}(r_\mend,\dend_2)\right|&\leq&
0.000252.
\end{eqnarray*}
Taking into account time rescaling $\rho=\ln r_\mbeg$ we
obtain the estimate
$$
\left|\frac{\partial}{\partial \rbeg}\tybeg(\rbeg,\dbeg_1)\right|=
\left|\frac{\partial}{\partial \rbeg}\unman(\ln \rbeg,\dbeg_1)\right| \leq \frac{7.65\cdot
10^{-6}}{r_1^-} =\frac{51}{21919}< 0.00233.
$$
Put $[\stman']:=0.000252\cdot[-1,1]$ and $[\tybeg']:=7.65\cdot 10^{-6}\cdot[-1,1]$.

Gathering this together we obtain a rigorous bound for the derivative
\begin{equation} \label{eq:proof-estm-F_1'}
\begin{array}{rcl}
F_1'([r_1^-,r_1^+])
		&\subset& \frac{\partial \pixend P^2_+}{\partial \rbeg}(W) + \frac{\partial \pixend P^2_+}{\partial \ybeg}(W)[\tybeg']\\
		&-&[\stman']\left(\frac{\partial \pi_{r} P^2_+}{\partial \rbeg}(W)+ \frac{\partial \pi_{r} P^2_+}{\partial \ybeg}(W)[\tybeg']\right)\\
		&\subset& [1538.59, 16826.7]+[-0.253485, 2.36953][\tybeg']\\
		&-& [\stman']\left([-292366, 565586]+ [-58.3544, 88.8462][\tybeg']\right)\\
		&\subset& [1396,16970],
\end{array}
\end{equation}
which is nonzero (see Remark~\ref{rem:estm-F_1'} for short comments about the possibility for obtaining  more precise bounds).

In a similar way we computed this derivative for $n=2,\ldots,6$ and we got the
following bounds
\begin{equation}\label{eq:rig-der-bounds}
\begin{array}{rcl}
F_2'([r_2^-,r_2^+]) & \subset & [-56203.7, -12850.2], \\
F_3'([r_3^-,r_3^+]) & \subset & [16342.1, 137481],\\
F_4'([r_4^-,r_4^+]) & \subset & [-183339, -24681.3],\\
F_5'([r_5^-,r_5^+]) & \subset & [39626.1, 307922], \\
F_6'([r_6^-,r_6^+]) & \subset & [-435119, -1138.93].
\end{array}
\end{equation}

It remains to show the transversality of $A_n$.

Observe that till now we established the following facts

\begin{itemize}
\item Let
$W^u_{(\rho \leq \rho_*),\Nbeg} \subset (-\infty,\rho_*] \times \Nbeg$ be a set consisting of all points $(\rho_0,z_0)$, such that there exists a solution $z(\rho)$ for $\rho\leq \rho_0$  of equation (\ref{eq:rho-neg}) (in coordinates in which $\Nbeg$ is defined this is the system (\ref{eq:beg-ord1})),
such that $z(\rho_0)=z_0$, $z((-\infty,\rho_0] \subset \Nbeg$ and $\lim_{\rho \to -\infty}z(\rho)=(0,0)$.

We have proved that there exists a smooth function $\unman(\rho,\xbeg)$, such that
\begin{equation*}
 W^u_{(\rho \leq \rho_*), \Nbeg}=\{(\rho,\xbeg, \unman(\rho,\xbeg)), \ | \  \rho \in (-\infty,\rho_*], \xbeg \in [-\dbeg_1,\dbeg_1]\}
 \end{equation*}

\item Let
$W^s_{(r \geq r_*),\Nend} \subset [r_*,\infty) \times \Nend$ be a set consisting of all points $(r_0,z_0)$, such that there exists a solution $z(r)$ for $r\geq r_0$  of equation (\ref{eq:end-ord1}) (in coordinates in which $\Nend$ is defined this is the system (\ref{eq:end-diag})),
such that $z(r_0)=z_0$, $z[r_0,\infty) \subset N_\mend$ and $\lim_{r \to \infty} z(r)=(0,0)$.

We have proved that there exists a smooth function $\stman(r,y_\mend)$, such that
\begin{equation*}
 W^s_{(r \geq r_*),\Nend}=\{(r,\stman(r,\yend),\yend), \ | \  r \in [r_*,\infty), y_\mend \in [-\dend_2,\dend_2]\}.
 \end{equation*}

\end{itemize}
The sets $ W^u_{(\rho \leq \rho_*),\Nbeg}$ and   $W^s_{(r \geq r_*),\Nend}$ are two dimensional surfaces in the extended phase space $\mathbb{R} \times \mathbb{R}^2$ (or $[0,\infty) \times \mathbb{R}^2$ depending whether we use the variable $\rho$ or $r$). These surfaces can be 'globalized' to invariant sets
$W^u$ and $W^s$,  by applying to $ W^u_{(\rho \leq \rho_*),\Nbeg}$ and   $W^s_{(r \geq r_*),\Nend}$ the flow induced by (\ref{eq:MS}). In this way we obtain two-dimensional immersed invariant manifolds.

Each of the solutions $A_n$ whose existence we already  established  belongs to $W^u \cap W^s$. We want to prove that the intersection along such solution is transversal.

Let $\gamma(r):=(r,A_n(r),A_n'(r))$  for some $n=1,\ldots,6$ be one of our connecting orbits. Let $p_\mbeg=(\barrbeg,\dbeg_1,\tybeg(\barrbeg,\dbeg_1))$ and  $p_\mend=(\barrend,\stman(\barrend,\dend_2),\dend_2)$ be the points on $\gamma$ on sections $\Pibeg$ (i.e. $\xbeg=\dbeg_1$) and $\Piend$ (i.e. $\yend=\dend_2$), respectively.  To prove the transversality
it is enough to show that the tangent spaces to $ W^u$ and $W^s$ at $p_\mend$ satisfy
\begin{equation*}
T_{p_\mend} W^u + T_{p_\mend} W^s = \mathbb{R}^3.
\end{equation*}
Since the point $p_\mend$ is in the isolating block $\Nend$ in coordinates $(\xend,\yend)$ we have
\begin{equation*}
  T_{p_\mend} W^s= T_{p_\mend} W^s_{(r \geq r_*),\Nend}= \mathrm{span} \left\{ \gamma'(\barrend), \left(1,\frac{\partial \stman}{\partial r}(\barrend,\dend_2) ,0 \right) \right\}
\end{equation*}
and these vectors are linearly independent because the last coordinate of $\gamma'(\barrend)$ is nonzero.
In the above formula the vector $\gamma'(\barrend)$ is the direction of the vector field (and the curve $\gamma$). The second vector
\begin{equation*}
  g_1=\left(1,\frac{\partial \stman}{\partial r}(\barrend,\dend_2) ,0 \right)
 \end{equation*}
  spans $T_{p_\mend} \left(W^s_{(r \geq r_*),\Nend} \cap \Piend \right)$, which is one dimensional.

 $T_{p_\mend} W^u$ is obtained from $T_{p_\mbeg} W^u_{(\rho \leq \rho^*),\Nbeg}$, as follows. $T_{p_\mbeg} W^u_{(\rho \leq \rho^*),\Nbeg}$ is spanned by $\gamma'(\barrbeg)$  and
 \begin{equation*}
 g_2=\left(1,0,\frac{\partial \tybeg }{\partial r} (\barrbeg,\dbeg_1) \right),
  \end{equation*}
  which spans the one-dimensional space $T_{p_\mbeg} \left(W^u_{(\rho \leq \rho_*),\Nbeg} \cap \Pibeg \right)$.

 Put $\tau=\barrend - \barrbeg$ which is the transition time for $P^{n+1}$ at $p_\mend$. The shift along the trajectory of $\gamma$ by $\tau$ maps $\gamma'(\barrbeg)$ to $\gamma'(\barrend)$. The Poincar\'e map $P^{n+1}$ acts on $g_2$ to produce a vector $\tilde{g}_2$ as follows
 \begin{eqnarray*}
    \pi_r \tilde{g}_2&=& \frac{\partial \pi_r P^{n+1}}{\partial r}(\barrbeg,\tybeg(\barrbeg,\dbeg_1)) + \frac{\partial \pi_r P^{n+1}}{\partial \ybeg}(\barrbeg,\tybeg(\barrbeg,\dbeg_1)) \frac{\partial \tybeg}{\partial r}(\barrbeg,\dbeg_1),  \\
    \pixend \tilde{g}_2&=& \frac{\partial \pixend P^{n+1}}{\partial r}(\barrbeg,\tybeg(\barrbeg,\dbeg_1)) + \frac{\partial \pixend P^{n+1}}{\partial \ybeg}(\barrbeg,\tybeg(\barrbeg,\dbeg_1)) \frac{\partial \tybeg}{\partial r}(\barrbeg,\dbeg_1), \\
    \piyend \tilde{g}_2&=&0.
 \end{eqnarray*}
 We have
 \begin{eqnarray*}
    T_{p_\mend} W^u= \mathrm{span} \left\{ \gamma'(\barrend), \tilde{g}_2   \right\}.
 \end{eqnarray*}
Therefore we need to show that the vectors $\gamma'(\barrend), g_1, \tilde{g}_2  $ are linearly independent.  Observe that $\piyend\gamma'(\barrend) \neq 0$, because the section $\Piend$ is transversal, and $\piyend g_1=\piyend \tilde{g}_2=0$. Therefore,
it is enough to check that the projections of $g_1,\tilde{g_2}$ on $(r,\xend)$ subspace are linearly independent.  This means the following
determinant is nonzero (we dropped arguments in the partial derivatives)
\begin{equation*}
D=\det \left[
  \begin{array}{cc}
    1, & \frac{\partial \pi_r P^{n+1}}{\partial r} + \frac{\partial \pi_r P^{n+1}}{\partial y} \frac{\partial \tybeg}{\partial r} \\
    \frac{\partial\xend}{\partial r}, & \frac{\partial \pixend P^{n+1}}{\partial r} + \frac{\partial \pixend P^{n+1}}{\partial y} \frac{\partial \tybeg}{\partial r}
  \end{array}
\right].
\end{equation*}
Observe that this leads to the following condition
\begin{equation*}
   D= \frac{\partial \pixend P^{n+1}}{\partial r} + \frac{\partial \pixend P^{n+1}}{\partial y} \frac{\partial \tybeg}{\partial r} -
    \frac{\partial \xend}{\partial r}  \left( \frac{\partial \pi_r P^{n+1}}{\partial r} + \frac{\partial \pi_r P^{n+1}}{\partial y} \frac{\partial \tybeg}{\partial r} \right)
    \neq 0.
\end{equation*}
We have
\begin{eqnarray*}
  \frac{\partial \pixend P^{n+1}}{\partial r} + \frac{\partial \pixend P^{n+1}}{\partial y} \frac{\partial \tybeg}{\partial r}=
  \frac{d}{dr}  \pixend P^{n+1}(r,\tybeg(r,\dbeg_2)), \\
  \frac{\partial \pi_r P^{n+1}}{\partial r} + \frac{\partial \pi_r P^{n+1}}{\partial y} \frac{\partial \tybeg}{\partial r} = \frac{d}{dr}
     \pi_{r} P^{n+1}(r,\tybeg(r,\dbeg_2)).
\end{eqnarray*}
therefore
\begin{equation*}
  D= F_n'(\bar{r}_\mend) \neq 0
\end{equation*}
hence the intersection is transversal.
\qed

\begin{rem}
\label{rem:estm-F_1'}
The estimation of $F_1'$ given by (\ref{eq:proof-estm-F_1'}) is very rough, but it is sufficient for our purposes. It is easy to see that the main reason for such big diameter of our estimate for $F_1'$ comes from our bound for $\frac{\partial \pi_{\xend} P^2}{\partial \rbeg}$. This number can be non-rigorously estimated by taking finite differences using data from (\ref{eq:cover}) (first row) and $\Delta r_1$ from (\ref{eq:r0-data}) to obtain 
$$ \frac{\pi_{\xend}P^2(r_1^+,\Dbeg_2) - \pi_{\xend}P^2(x_1^-,\Dbeg_2)}{2 \Delta r_1}\approx (0.003294+0.003223)/(8 \cdot 10^{-7})=8146.25. $$
Subdividing $\Delta r_1$ onto $100$ pieces we got a rigorous estimate on the derivative
$$F_1'([r_1^-,r_1^+])\subset[7902.54, 8399.24].$$

From the above computations it is clear that the key ingredients in obtaining nonzero derivative $F'([r_1^-,r_1^+])$ are
  sharp bounds for $[\xend']$ obtained from the cone condition and the transversality coming from dynamics, i.e. large in absolute value coefficient $\frac{\partial \pixend P_+^2}{\partial \rbeg}(W)$.
\end{rem}

\subsection{Implementation notes.}
A short C++11 program that realizes part of a computer assisted proof of Theorem~\ref{thm:six-orbits} is available from \cite{cpp}.
All the inequalities, i.e. return times (\ref{eq:return-time}), shooting (\ref{eq:cover}) and estimation on derivatives (\ref{eq:rig-der-bounds}), regarding all six connecting orbits are checked within 0.7 second on a laptop type computer with Intel Core i7 2GHz processor. The program has been compiled and tested under Ubuntu OS with gcc-4.8.1 and gcc-4.9.1, OSX 10 with clang 6.0 and MS Windows 7 with gcc-4.8.1.

This program does not contain nonrigorous routines for finding good candidates for connecting orbits. They have been found by a simple  bisection algorithm. The sizes of sets (see (\ref{eq:r0-data})) were adjusted by hand so that no subdivisions are necessary when computing derivatives $F_n'$.

\end{document}